\documentclass[12pt]{article}
\usepackage[left=20mm, top=20mm, right=20mm, bottom=20mm]{geometry}
\usepackage{amsmath,amsthm, amssymb}
\usepackage[font=scriptsize]{caption}
\usepackage{multicol}
\usepackage{graphicx}
\usepackage{footmisc}
\usepackage{float}
\usepackage{geometry}
\usepackage{longtable}
\usepackage{makeidx}
\usepackage{endnotes}
\usepackage[norefeq]{nomencl}
\usepackage[latin9]{inputenc}
\usepackage{esint}
\usepackage[justification=centering]{caption}
\makeglossary

\begin{document}

\title{\textbf{An Improved Non-linear Weights for Seventh-Order WENO Scheme}}
\author{Samala Rathan\thanks{Email: rathan.maths@gmail.com},\hspace{1mm}
G Naga Raju \thanks{Email: gnagaraju@mth.vnit.ac.in}\\
 {\small{}{}{}Department of Mathematics, Visvesvaraya National Institute
of Technology, Nagpur, India}}
\date{}
\maketitle
\begin{abstract}
In this article, the construction and implementation of a seventh
order weighted essentially non-oscillatory scheme is reported for
hyperbolic conservation laws. Local smoothness indicators are constructed
based on $L_{1}$-norm, where a higher order interpolation polynomial
is used with each derivative being approximated to the fourth order
of accuracy with respect to the evaluation point. The global smoothness
indicator so constructed ensures the scheme achieves the desired order
of accuracy. The scheme is reviewed in the presence of critical points
and verified the numerical accuracy, convergence with the help of
linear scalar test cases. Further, the scheme is implemented to non-linear
scalar and system of equations in one and two dimensions. As the formulation
is based on method of lines, to move forward in time linear strong-stability-preserving
Runge-Kutta scheme $\left(lSSPRK\right)$ for the linear problems
and the fourth order nonlinear version of five stage strong stability
preserving Runge-Kutta scheme $\left(SSPRK(5,4)\right)$ for nonlinear
problems is used.
\end{abstract}

\section{Introduction}

\hspace{0.6 cm}The hyperbolic conservation laws arise in many applications
such as in gas dynamics, magnetohydrodynamics (MHD) and shallow water
flows. It is well known that even if the initial conditions are smooth,
the hyperbolic conservation laws may develop discontinuities in its
solution, such as shocks, contact discontinuities etc. Godunov \cite{Godunov}
was first to propose a first order upwind scheme for the solution
of these equations in the year $1959,$ which turned out to be a stepping
stone for the development of various upwind schemes in the following years.
In order to construct a higher order scheme Harten \cite{Ha83,Ha84}
introduced the concept of Total Variation Diminishing (TVD), which
says that the total variation to the approximation of numerical solution
must be non-increasing with time. Later it was shown that the TVD
schemes are having at most of first-order accuracy near smooth extrema
\cite{OC84} .

Harten et al. \cite{HOEC86,HO87,HO97} derived the higher order schemes
with the property of relaxing the TVD condition and allowing the occurrence
of spurious oscillations in the numerical scheme but the $O(1)$ Gibbs-like
phenomena is essentially prevented, which is termed as Essentially
Non-oscillatory (ENO) property and the schemes are known as ENO schemes.
These are the first successful higher order schemes for the spatial
discretization of the hyperbolic conservation laws, in a finite-volume
formulation. The ENO schemes adopts a strategy of choosing the interpolation
points over a stencil which avoids the induction of oscillations in
the numerical solution through a smoothness indicator of a solution.
And based on this idea the smoothest stencil is chosen from a set
of candidate stencils. As a result, the ENO scheme obtains information
from smooth regions and avoids spurious oscillations near discontinuities.
Further, these schemes were studied in a finite-difference environment
by Shu and Osher \cite{osherSHU,Shu-osher1}.

The weighted ENO (WENO) schemes are set forth by Liu et al. \cite{XD Liu},
in a finite-volume frame of reference up to third-order of accuracy.
Later, Jiang and Shu \cite{Jiang and shu7} have put forward these
WENO schemes in a finite-difference setup to a higher order accuracy
with the new smoothness indicators. These smoothness indicators are
measured in the scaled $L_{2}$ norm, that is, they are the sum of
the normalized squares, of all derivatives of the local interpolating
polynomials. This scheme is referred as WENO-JS in the content to
follow. For more details on ENO and WENO schemes, one can refer to
the articles \cite{shu notes} and \cite{c w shu16}. A very high
order schemes are constructed in a similar manner of WENO-JS in \cite{balsara and shu8},
which we mention them here as WENO-BS schemes. Seventh order WENO-BS
scheme is revised in \cite{ShenZha,ShenZha2} and inspected the scheme
in the presence of critical points.

Henrick et al. \cite{henrick aslam powers5} examined that the actual
convergence rate of the fifth-order WENO-JS scheme is less than the
desired order, for the problems where the first and third order derivatives
of the flux do not vanish simultaneously. In addition, it was ascertained
that the convergence rate of the scheme is sensitive to the parameter
$\epsilon,$ employed in the evaluation of smoothness indicators to
overcome from vanishing denominator. The authors revived the WENO-JS
scheme by using a mapping function on the nonlinear weights such that
the scheme, named as mapped WENO, satisfies the sufficient condition
where WENO-JS fails and achieves an optimal order of convergence near
simple smooth extrema. Subsequently, a very high order WENO schemes
were developed based on the mapping function in \cite{Gerolymus}.

Borges et al. \cite{Borges caramona10} reviewed the fifth-order WENO
schemes, entitled as WENO-Z scheme, by initiating a global smoothness
indicator, which measures the smoothness of the larger stencil utilized
in the construction of nonlinear weights. It was numerically validated
that WENO-Z scheme is less dissipative than WENO-JS scheme and more
efficient than mapped WENO scheme. WENO-Z scheme retained the convergence
order as four at the first-order critical points, degrade to two when
higher order critical points are encountered. These thoughts are extended
by Castro et al. \cite{Castro et al} to higher order schemes and
produced a closed-form formula for the global smoothness indicators.
The authors also assessed the dominance of the parameters $p$ and
$\epsilon$ to retain the desired order of accuracy. The parameter
$p$ is set up in the formulation of nonlinear weights to ascertain
that these nonlinear weights converge to the ideal weights at a fast
enough convergence rate.

The convergence analysis of WENO-JS scheme explored by Arandiga et
al. \cite{Arandiga} is based on the value of $\epsilon,$ proposed
that $\epsilon$ value is proportional to the square of mesh size
$\Delta x$, instead of a constant value so that the scheme achieves
$(2r-1)^{th}$ order of accuracy at smooth regions regardless of neighboring
extrema, while this is of order $r$ when the function has a discontinuity
in the stencil of $(2r-1)$ points and is smooth in at least one of
the $r-$point stencil. A question about the behavior of WENO-Z scheme
when the $\epsilon$ value is taken in accordance with the value mentioned
in \cite{Arandiga}, is examined by Don and Borges \cite{W S DON }.
The authors made the accuracy analysis of the WENO-Z scheme and suggested
a condition on the value of $\epsilon$ to achieve the full global-order
of accuracy as similar to that of \cite{Arandiga}. Further the authors
have shown that the numerical oscillations can be attenuated by increasing
the parameter $p$ value from $2$ to $r-1$.

An alternate to the smoothness indicators of fifth order WENO-JS scheme
were formulated by Fan et al. \cite{P fan et.al.14} with the help
of Lagrange interpolation polynomials, accordingly a very high order
schemes were derived by Fan in \cite{P fan15}. These schemes are
based on the idea of constructing higher order global smoothness indicators,
due to which less dissipation occurs in the solution near discontinuities.
Very recently, another version of a smoothness indicators were proposed
by Ha et al. \cite{Ha et.al13}, measured in $L_{1}$-norm, hereafter
referred as WENO-NS scheme. The authors introduced a higher-order
approximation to the first derivative in the formation of local smoothness
indicators which yields an improved behavior relative to other fifth-order
WENO schemes. The global smoothness indicator for the WENO-NS scheme
is preferred as an average of the two measurements, the smoothness
information of the five point stencil and the middle three point stencil.

Kim et al. \cite{KIM et al} perceived that, the three sub-stencils
of the fifth-order WENO-NS scheme provides an unbalanced contribution
to the flux at an evaluation point along the interface and an additional
contribution term which measures the smoothness of the middle stencil
in the formation of global smoothness indicator. The authors made
a balanced tradeoff among the sub-stencils through a parameter $\delta$
and a global smoothness indicator is figured out which doesn't depend
on the smoothness information of the middle three point stencil anymore.
These modifications lead to better results than the WENO-NS as well
as other fifth-order WENO schemes, this scheme is pointed out as WENO-P
scheme in the later part of this article.

A simple analysis verifies that the WENO-NS and WENO-P schemes attain
the fifth-order accuracy at first order critical points but fails
to achieve the accuracy at the points where second derivative vanishes.
We have suggested a modified WENO-P scheme in \cite{RathanRaju} based
on the idea of the linear combination of second-order derivatives,
leading to a higher-order derivative information, is used in the construction
of a global smoothness indicator. The modified smoothness indicator
satisfies the sufficient condition, assert the requirement to achieve
desired order of accuracy, even in the case of second order derivative
vanishes.

In this article, a seventh order WENO scheme is derived in the lines
of \cite{Ha et.al13} and \cite{RathanRaju}. The smoothness indicators
are obtained from the generalized undivided difference operator. Each
of this operator is up to fourth-order of accuracy at the evaluation
point, so the resulting scheme is seventh order accurate. Introduced
parameters $\xi_{1}$, $\xi_{2}$ to balance the tradeoff between
the accuracies around the smooth to the discontinuous regions. The
global smoothness indicator so earned satisfies the sufficient condition
to get the optimal order of convergence rate, unvarying in the presence
of critical points. Utilized strong stability preserving Runge-Kutta
schemes introduced in \cite{Gottilieb} to advance the time. These
are detailed out in the following sections, briefly they are:

Section $2$, deals with the preliminaries about WENO reconstruction
to the one-dimensional scalar conservation laws and section $3$ introduces
the proposed WENO scheme where a new global measurement and local
smoothness indicators are derived, which estimates the smoothness
of a local solution in the construction of a seventh-order WENO scheme.
Numerics for one-dimensional scalar test problems such as linear advection,
Burger's equation, the examples pertaining to the system of Euler
equations such as shock tube problems, 1D shock-entropy wave interaction
problem and 2D Riemann problem of gas dynamics are reported in Section
$4$, to demonstrate the advantages of the proposed WENO scheme. Finally,
concluding remarks are in Section $5$.

\section{The fundamentals of WENO scheme}

\hspace{0.6 cm}This part accounts to the construction of the seventh
order weighted essentially non-oscillatory scheme in a finite difference
framework for approximate the solution of hyperbolic conservation
laws
\begin{equation}
u_{t}+f(u)_{x}=0,\,x\in(-\infty,\infty),\,t\in(0,\infty),\label{eq:1}
\end{equation}
with initial condition
\begin{equation}
u(x,0)=u_{0}(x),\,x\in(-\infty,\infty).\label{eq:IC}
\end{equation}
Here $u=(u_{1},u_{2},......,u_{m})^{T}$ is a $m-$dimensional vector
of conserved variables defined for space and time variables $x\text{ and }t$
respectively, $f$ is a flux function which depends on the conserved
quantity $u.$ The system $(\ref{eq:1})$ is called hyperbolic if
all the eigen values $\lambda_{1},\lambda_{2},...,\lambda_{m}$ of
the Jacobian matrix $A=\frac{\partial f}{\partial u}$ of the flux
function are real and the set of right eigen vectors are complete.

For numerical approximation the spatial domain is discretized with
uniform grid, for brevity in the presentation. Let $\triangle x=x_{j+\frac{1}{2}}-x_{j-\frac{1}{2}}$
be the length of the $j^{\text{th}}$ cell $I_{j}=\left[x_{j-\frac{1}{2}},x_{j+\frac{1}{2}}\right]$
with center $x_{j}=\frac{1}{2}\left(x_{j+\frac{1}{2}}+x_{j-\frac{1}{2}}\right),$
here $x_{j\pm\frac{1}{2}}$ are known as cell interfaces. The approximation
of the spatial derivative in the hyperbolic conservation laws $(\ref{eq:1})$
yields a semi-discrete formulation
\begin{equation}
\frac{d{u_{j}}}{dt}=-\frac{1}{\triangle x}\left({\hat{f}_{j+\frac{1}{2}}-\hat{f}_{j-\frac{1}{2}}}\right).\label{eq:2}
\end{equation}
Here ${u_{j}}$ is an approximation to $u$ at a point $x=x_{j}$
in time $t,$ i.e., for the value $u(x_{j},t)$ and $\hat{f}_{j\pm\frac{1}{2}}$
are the numerical fluxes which are Lipschitz continuous in each of
its arguments and are consistent with the physical flux, that means
$\hat{f}(u,....,u)=f(u).$ The conservation property is retrieved
by defining a function $h$ implicitly through the equation (see Lemma
$2.1$ of \cite{Shu-osher1})
\begin{equation}
f(x):= f\left(u(x,.)\right)=\frac{1}{\triangle x}\intop_{x-\frac{\triangle x}{2}}^{x+\frac{\triangle x}{2}}h(\xi)d\xi.\label{eq:3}
\end{equation}
Differentiating (\ref{eq:3}) with respect to $x$ yields
\[
f\left(u(x,.)\right)_{x}=\frac{1}{\Delta x}\left(h\left(x+\frac{\Delta x}{2}\right)-h\left(x-\frac{\Delta x}{2}\right)\right),
\]
thus $h\left(x\pm\frac{\Delta x}{2}\right)$ should be an approximation
to the numerical flux $\hat{f}_{j\pm\frac{1}{2}},$ such that
\[
\hat{f}_{j\pm\frac{1}{2}}=h\left(x_{j\pm\frac{1}{2}}\right)+O\left(\Delta x^{2r-1}\right),
\]
$r$ represents the number of cells in a stencil. Thus the numerical
flux $\hat{f}_{j\pm\frac{1}{2}}$ can be acquired by using higher
order polynomial interpolation to $h\left(x_{j\pm\frac{1}{2}}\right)$
with the help of known values of $f(x)$ at the cell centers, $f_{j}=f\left(x_{j}\right).$

To ensure the numerical stability and to avoid entropy violating solutions,
the flux $f(u)$ is splitted into two parts $f^{+}$ and $f^{-},$
the positive and negative parts of $f(u)$ respectively, such that
\begin{equation}
f(u)=f^{+}(u)+f^{-}(u),\label{eq:3-1}
\end{equation}
where $\frac{df^{+}(u)}{du}\geq0$ and $\frac{df^{-}(u)}{du}\leq0.$
The numerical fluxes $\hat{f}^{+}$ and $\hat{f}^{-}$ evaluated at
$x=x_{j+\frac{1}{2}}$ reduces (\ref{eq:3-1}) as
\[
\hat{f}_{j+\frac{1}{2}}=\hat{f}_{j+\frac{1}{2}}^{+}+\hat{f}_{j+\frac{1}{2}}^{-}.
\]
We will describe here how $\hat{f}_{j+\frac{1}{2}}^{+}$ can be approximated,
as $\hat{f}_{j+\frac{1}{2}}^{-}$ is symmetric to the positive part
with respect to $x_{j+\frac{1}{2}}.$ In the description for the approximation
of $\hat{f}_{j+\frac{1}{2}}^{+}$ to follow we drop the $^{,}+^{,}$
sign in the superscript, for simplicity.

\subsection{Seventh order WENO scheme}

\hspace{0.6 cm}WENO scheme prefers $(2r-1)$ points global stencil,
$S^{2r-1},$ to achieve $(2r-1)^{th}$ order of accuracy. The stencil
$S^{2r-1}$ is subdivided into $r$ sub-stencils with each sub-stencil
bearing $r$ cells. In particular, seventh-order WENO scheme accounts
to a $7-$points stencil, which is subdivided into four $4$-points
sub-stencils $S^{0},S^{1},S^{2},S^{3}.$ In accordance with cell $I_{j},$
each sub-stencil encloses four grid points, specified as
\[
S^{k}(j)=\{x_{j+k-3},x_{j+k-2},x_{j+k-1},x_{j+k}\},\,k=0,1,2,3.
\]
A third degree interpolating polynomial $\hat{f}^{k}(x)$ is formulated
in each sub-stencil $S^{k}(j)$ and evaluating it at the cell boundary
$x_{j+\frac{1}{2}},$ we retain
\begin{equation}
\hat{f}^{k}\left(x_{j+\frac{1}{2}}\right)=\hat{f}_{j+\frac{1}{2}}^{k}=\sum_{q=0}^{3}c_{k,q}f_{j+k-3+q},\label{f^kjph}
\end{equation}
where the coefficients $c_{k,q}(q=0,1,2,3)$ are the Lagrange interpolation
coefficients, independent of the values of the flux function $f,$
but depends on the left shift parameter $k=0,1,2,3.$ The equation
(\ref{f^kjph}) on each stencil takes the form
\begin{align*}
 & \hat{f}_{j+\frac{1}{2}}^{0}=-\frac{1}{4}f_{j-3}+\frac{13}{12}f_{j-2}-\frac{23}{12}f_{j-1}+\frac{25}{12}f_{j},\\
 & \hat{f}_{j+\frac{1}{2}}^{1}=\frac{1}{12}f_{j-2}-\frac{5}{12}f_{j-1}+\frac{13}{12}f_{j}+\frac{1}{4}f_{j+1},\\
 & \hat{f}_{j+\frac{1}{2}}^{2}=-\frac{1}{12}f_{j-1}+\frac{7}{12}f_{j}+\frac{7}{12}f_{j+1}-\frac{1}{12}f_{j+2},\\
 & \hat{f}_{j+\frac{1}{2}}^{3}=\frac{1}{4}f_{j}+\frac{13}{12}f_{j+1}-\frac{5}{12}f_{j+2}+\frac{1}{12}f_{j+3}.
\end{align*}
The fluxes $\hat{f}_{j-\frac{1}{2}}^{k}$ can be fetched through shifting
the index to the left by one in (\ref{f^kjph}). The Taylor's expansion
of the fluxes $\hat{f}_{j\pm\frac{1}{2}}^{k},$ settle in as
\[
\hat{f}_{j\pm\frac{1}{2}}^{k}=h_{j\pm\frac{1}{2}}+A_{4}^{k}f_{j}^{(4)}\Delta x^{4}+O\left(\Delta x^{5}\right),
\]
where $A_{4}^{k}$ is the leading order coefficient in the expansion
and $f_{j}^{(4)}$ is the $4^{th}$ derivative of $f(x)$ at $x=x_{j}.$
The convex combination of these flux functions leads to the approximation
of $\hat{f}_{j+\frac{1}{2}},$ that is, we set
\begin{equation}
\hat{f}_{j+\frac{1}{2}}=\sum_{k=0}^{3}\omega_{k}\hat{f}_{j+\frac{1}{2}}^{k}.\label{eq:recon}
\end{equation}
Here $\omega_{k}^{,}s$ are non-linear weights, satisfying the conditions

\[
\omega_{k}\geq0,\,k=0,1,2,3\text{ and }\,\sum_{k=0}^{3}\omega_{k}=1.
\]
If the function $h(x)$ is free from discontinuities in all of the
sub-stencils $S^{k}(j),$ $k=0,1,2,3,$ we can assess the constants
$d_{k}$ such that the linear combination of $\hat{f}_{j\pm\frac{1}{2}}^{k}$
provides the seventh order convergence to $\hat{f}_{j\pm\frac{1}{2}},$
that is,
\[
\hat{f}_{j\pm\frac{1}{2}}=\sum_{k=0}^{3}d_{k}\hat{f}_{j\pm\frac{1}{2}}^{k}=h_{j\pm\frac{1}{2}}+A_{7}f_{j}^{(7)}\Delta x^{7}+O\left(\triangle x^{8}\right).
\]
The $d_{k}^{,}s$ are termed as the ideal weights since they invokes
the upstream central scheme of seventh-order, in seven points stencil.
The values of these ideal weights are evaluated as
\begin{equation}
d_{0}=1/35,\,d_{1}=12/35,\,d_{2}=18/35,\,d_{3}=4/35.\label{eq:idealwts}
\end{equation}
When the non-linear weights $\omega_{k}$ are equal to the ideal weights
$d_{k},$ we have
\[
\frac{1}{\Delta x}\left(\hat{f}_{j+\frac{1}{2}}-\hat{f}_{j-\frac{1}{2}}\right)=\frac{1}{\Delta x}\left(\left(h_{j+\frac{1}{2}}-h_{j-\frac{1}{2}}\right)+O\left(\triangle x^{8}\right)\right)=f^{\prime}\left(x_{j}\right)+O\left(\triangle x^{7}\right),
\]
thus the approximation to the spatial derivative of the flux $f(x)$
at the cell-centre $x=x_{j}$ is $O\left(\triangle x^{7}\right).$
Hence the sufficient condition to achieve the seventh order convergence
for the scheme is given by
\begin{equation}
\omega_{k}^{\simeq}-d_{k}=O\left(\triangle x^{4}\right)\label{eq:17}
\end{equation}
where the superscripts $^{,}\sim{}^{,}$ and $^{,}-{}^{,}$  on $\omega_{k}$
corresponds to their use in either $\hat{f}_{j+\frac{1}{2}}^{k}$
and $\hat{f}_{j-\frac{1}{2}}^{k}$ respectively. Note that the condition
$(\ref{eq:17})$ can be weakened for WENO-JS scheme to $\omega_{k}^{\simeq}-{d}_{k}=g\left(x_{j+\frac{1}{2}}\right)\left(\triangle x^{3}\right)+O\left(\triangle x^{4}\right)$
for a locally Lipschitz continuous function $g$ under a suitable
condition, which depends on the value of $\epsilon.$

\section{The numerical scheme, WENO-NS7}

\hspace{0.6 cm}The essential ingredient of the WENO schemes is in
the computation of smoothness indicators, Ha et. al. \cite{Ha et.al13}
have established the smoothness indicators measured in $L_{1}-$norm
for the fifth-order WENO scheme. Here we are extending it to the higher
order schemes, particularly for the seventh order WENO scheme and
calling the scheme as WENO-NS7. The primary thought is of getting
a higher order approximation to the derivatives, as it is known that
a smoothness indicators based on $L_{1}$-norm may give a loss of
accuracy in smooth regions \cite{serna}. The heuristic construction
of smoothness indicators is as follows.

Define the operators $L_{s,k}f$ which measures the regularity of
the solution in each of the four points stencil $S^{k}(j)=\{x_{j+k-3},x_{j+k-2},x_{j+k-1},x_{j+k}\},\,k=0,1,2,3$
by estimating the approximate magnitude of derivatives. Once obtained
these operators, the smoothness indicators $\beta_{k}$ are designed
as
\begin{equation}
\beta_{k}=\xi_{1}\left|L_{1,k}f\right|+\xi_{2}\left|L_{2,k}f\right|+\left|L_{3,k}f\right|,\label{eq:beta_k}
\end{equation}
here $\:\xi_{i}\in(0,1],\:i=1,2,$ are the parameters introduced to
balance the tradeoff between the accuracy around the smooth region
to that of the discontinuous region. The operators $L_{s,k}f$ are
the generalized undivided differences of $f,$ which are formulated
as
\begin{equation}
L_{s,k}f\left(x_{j+\frac{1}{2}}\right):=\sum_{x_{l}\in S^{k}(j)}c_{k,l}^{[s]}f\left(x_{l}\right),\,s=1,2,3,\label{eq:opL}
\end{equation}
where the coefficient vector
\begin{equation}
c_{k}^{[s]}:=\left(c_{k,l}^{[s]}:x_{l}\in S^{k}(j)\right)^{T},\,s=1,2,3,\label{26}
\end{equation}
is obtained by solving the linear system
\[
\sum_{x_{l}\in S_{k}(j)}c_{k,l}^{[s]}\frac{\left(x_{l}-x_{j+\frac{1}{2}}\right)^{m}}{m!}=\begin{cases}
\Delta x^{s}, & \text{ if }s=m,\\
0, & \text{ if }s\neq m.
\end{cases}
\]
with $m=0,1,2,3$. This linear system can be written in the matrix
form
\[
Mc_{k}^{[s]}=r,
\]
with the matrices $M$ and $r$ defined by
\begin{eqnarray*}
M= & M_{k,j}:= & \left(\frac{\left(x_{l}-x_{j+\frac{1}{2}}\right)^{m}}{m!}:x_{l}\in S_{k}(j),m=0,1,2,3\right),\\
r= & r_{s}:= & \left(\Delta x^{s}\delta_{s,m}:m=0,1,2,3\right)^{T},
\end{eqnarray*}
where $\delta_{s,m}$ is the Kronecker delta and there exist a unique
solution for this linear system, as $M$ is a Vandermode matrix. With
the coefficients (\ref{26}) the operators (\ref{eq:opL}) takes the
form
\begin{align*}
 & L_{1,0}f=\frac{1}{24}(-23f_{j-3}+93f_{j-2}-141f_{j-1}+71f_{j}),\,L_{2,0}f=\frac{1}{2}(-3f_{j-3}+11f_{j-2}-13f_{j-1}+5f_{j}),\\
 & L_{1,1}f=\frac{1}{24}(f_{j-2}-3f_{j-1}-21f_{j}+23f_{j+1}),\,L_{2,1}f=\frac{1}{2}(-3f_{j-2}+5f_{j-1}-7f_{j}+3f_{j+1}),\\
 & L_{1,2}f=\frac{1}{24}(f_{j-1}-27f_{j}+27f_{j+1}-f_{j+2}),\,L_{2,2}f=\frac{1}{2}(f_{j-1}-f_{j}-f_{j+1}+f_{j+2}),\\
 & L_{1,3}f=\frac{1}{24}(-23f_{j}+21f_{j+1}+3f_{j+2}-f_{j+3}),\,L_{2,3}f=\frac{1}{2}(3f_{j}-7f_{j+1}+5f_{j+2}-f_{j+3}),
\end{align*}
\begin{align*}
 & L_{3,0}f=(-f_{j-3}+3f_{j-2}-3f_{j-1}+f_{j}),\\
 & L_{3,1}f=(-f_{j-2}+3f_{j-1}-3f_{j}+f_{j+1}),\\
 & L_{3,2}f=(-f_{j-1}+3f_{j}-3f_{j+1}+f_{j+2}),\\
 & L_{3,3}f=(-f_{j}+3f_{j+1}-3f_{j+2}+f_{j+3}).
\end{align*}

The third operator $L_{3,k}f$ is the same as in the WENO-BS scheme
which is described in appendix. However WENO-NS7 scheme uses the absolute
values where as the WENO-BS uses the squared ones. The advantages
with these operators $L_{s,k}f$ is that the approximation of the
derivative $\Delta x^{s}f^{(s)}$ to be of higher order accuracy at
the evaluation point, $x_{j+\frac{1}{2}},$ which is stated in the
following theorem.\\
\textbf{Thereom}
Let the stencil $S^{k}(j)=\{x_{j+k-3},x_{j+k-2},x_{j+k-1},x_{j+k}\},\,k=0,1,2,3,$
and assume that $f\in C^{4}(I)$ where $I$ is an open interval containing
$S_{k}(j)$. For each $s=1,2,3,$ the operator (\ref{eq:opL}) satisfies
\[
L_{s,k}f(x_{j+\frac{1}{2}})=\frac{d^{s}f}{dx^{s}}\Delta x^{s}+O\big(\Delta x^{4}\big).
\]
\textbf{proof}
The proof follows in a similar lines that of the Theorem 3.2 of \cite{Ha et.al13}.
The Taylor's expansion of the operators $L_{s,k}f$ for each $s$
and $k$ reveals
\[
\begin{array}{cc}
 & L_{1,0}f=\Delta xf_{j+\frac{1}{2}}^{(1)}-\frac{22}{24}\Delta x^{4}f_{j+\frac{1}{2}}^{(4)}+O(\Delta x^{5}),\,L_{2,0}f=\Delta x^{2}f_{j+\frac{1}{2}}^{(2)}-\frac{43}{24}\Delta x^{4}f_{j+\frac{1}{2}}^{(4)}+O(\Delta x^{5}),\\
 & L_{1,1}f=\Delta xf_{j+\frac{1}{2}}^{(1)}+\frac{1}{24}\Delta x^{4}f_{j+\frac{1}{2}}^{(4)}+O(\Delta x^{5}),\,L_{2,1}f=\Delta x^{2}f_{j+\frac{1}{2}}^{(2)}-\frac{7}{24}\Delta x^{4}f_{j+\frac{1}{2}}^{(4)}+O(\Delta x^{5}),\\
 & L_{1,2}f=\Delta xf_{j+\frac{1}{2}}^{(1)}+O(\Delta x^{5}),\qquad\quad\qquad L_{2,2}f=\Delta x^{2}f_{j+\frac{1}{2}}^{(2)}+\frac{5}{24}\Delta x^{4}f_{j+\frac{1}{2}}^{(4)}+O(\Delta x^{5}),\\
 & L_{1,3}f=\Delta xf_{j+\frac{1}{2}}^{(1)}-\frac{1}{24}\Delta x^{4}f_{j+\frac{1}{2}}^{(4)}+O(\Delta x^{5}),\,L_{2,3}f=\Delta x^{2}f_{j+\frac{1}{2}}^{(2)}-\frac{7}{24}\Delta x^{4}f_{j+\frac{1}{2}}^{(4)}+O(\Delta x^{5}),
\end{array}
\]
\begin{align*}
 & L_{3,0}f=\Delta x^{3}f_{j+\frac{1}{2}}^{(3)}-2\Delta x^{4}f_{j+\frac{1}{2}}^{(4)}+O(\Delta x^{5}),\\
 & L_{3,1}f=\Delta x^{3}f_{j+\frac{1}{2}}^{(3)}-\Delta x^{4}f_{j+\frac{1}{2}}^{(4)}+O(\Delta x^{5}),\\
 & L_{3,2}f=\Delta x^{3}f_{j+\frac{1}{2}}^{(3)}+O(\Delta x^{5}),\\
 & L_{3,3}f=\Delta x^{3}f_{j+\frac{1}{2}}^{(3)}-\Delta x^{4}f_{j+\frac{1}{2}}^{(4)}+O(\Delta x^{5}).
\end{align*}
Thus the operators are in tuned with the Theorem 3.1 stated above.

Now the non-linear weights for the scheme are defined as
\begin{equation}
\omega_{k}^{NS7}=\frac{\alpha_{k}^{NS7}}{\sum_{q=0}^{3}\alpha_{q}^{NS7}},\,\,\,\alpha_{k}^{NS7}=d_{k}\bigg(1+\frac{\zeta}{(\beta_{k}+\epsilon)^{2}}\bigg),\,k=0,1,2,3,\label{eq:nlwts}
\end{equation}
where $\zeta,$ a global smoothness indicator is taken as
\begin{equation}
\zeta=\left|\beta_{0}-\beta_{3}\right|^{2},\label{eq:gs}
\end{equation}
and $d_{k}$ the ideal weights, given in (\ref{eq:idealwts}). To
avoid the scenario of zero division a small number $0<\epsilon\ll1$
is incorporated in the calculations of non-linear weights (\ref{eq:nlwts}).

Next we discuss the convergence analysis of the WENO-NS7 scheme, in-particularly
at the critical points, i.e., we analyze how the weights $\omega_{k}^{NS7}$
approaches to the ideal weights $d_{k}$ in the presence of critical
points.

\subsection{Convergence order of WENO-NS7 scheme}

\hspace{0.6 cm}First consider that there are no critical points and
let's take $\epsilon=0$ in (\ref{eq:beta_k}), from Taylor's series
expansion
\begin{equation}
{\beta_{k}}=\xi_{1}\left|\Delta xf_{j+\frac{1}{2}}^{(1)}\right|+\xi_{2}\left|\Delta x^{2}f_{j+\frac{1}{2}}^{(2)}\right|+\left|\Delta x^{3}f_{j+\frac{1}{2}}^{(3)}\right|+O\left(\Delta x^{4}\right)\text{ for all }k.\label{eq:beta_k exp}
\end{equation}
Similarly from the expansion of the global smoothness indicator (\ref{eq:gs}),
we've
\begin{equation}
\zeta=\left[\frac{81}{24}\Delta x^{4}f_{j+\frac{1}{2}}^{(4)}+O\left(\Delta x^{5}\right)\right]^{2}.\label{eq:gsexp}
\end{equation}
Then there exist a constant $D$ such that
\begin{align}
1+\frac{\zeta}{{{\beta_{k}}^{2}}} & =1+D\Delta x^{6}+O\left(\Delta x^{7}\right),\label{eq:1c}\\
 & =(1+D\Delta x^{6})\bigg(1+\frac{O(\Delta x^{7})}{1+D\Delta x^{6}}\bigg),\nonumber \\
 & =D_{\Delta x}\left(1+O\left(\Delta x^{7}\right)\right),\nonumber
\end{align}
where $D_{\Delta x}=(1+D\Delta x^{6})>0$.

If $f_{j+\frac{1}{2}}^{(1)}=0,\;f_{j+\frac{1}{2}}^{(2)}\neq0$, from
(\ref{eq:beta_k exp}) the smoothness indicators are of the form
\[
{\beta_{k}}=\left|\Delta x^{2}f_{j+\frac{1}{2}}^{(2)}\right|\left(1+O(\Delta x)\right),
\]
then there exists a constant $D_{1}$ such that
\begin{align}
1+\frac{\zeta}{{{\beta_{k}}^{2}}} & =1+D_{1}\Delta x^{4}+O(\Delta x^{5}),\label{eq:2c}\\
 & =(1+D_{1}\Delta x^{4})\bigg(1+\frac{O(\Delta x^{5})}{1+D_{1}\Delta x^{4}}\bigg),\nonumber \\
 & =D_{\Delta x}^{1}(1+O(\Delta x^{5})),\nonumber
\end{align}
where $D_{\Delta x}^{1}=(1+D_{1}\Delta x^{4})>0$. Similarly if $f_{j+\frac{1}{2}}^{(1)}=0,f_{j+\frac{1}{2}}^{(2)}=0,f_{j+\frac{1}{2}}^{(3)}\neq0,$
then there exist a constant $D_{2}$ such that
\begin{equation}
1+\frac{\zeta}{{{\beta_{k}}^{2}}}=1+D_{2}\Delta x^{2}+O(\Delta x^{3}).\label{eq:3c}
\end{equation}

From (\ref{eq:1c}) and (\ref{eq:2c}), the weights $\omega_{k}^{NS7}$
satisfies the sufficient condition (\ref{eq:17}) if the first derivative
vanishes but not the second derivative. In order to satisfy the sufficient
condition even at higher order critical points, the nonlinear weights
are defined as
\begin{equation}
\omega_{k}^{NS7}=\frac{\alpha_{k}^{NS7}}{\sum_{q=0}^{3}\alpha_{q}^{NS7}},\;\alpha_{k}^{NS7}=d_{k}\left[1+\left(\frac{\zeta}{\left(\beta_{k}+\epsilon\right)^{2}}\right)^{s}\right],\,k=0,1,2,3,\label{eq:3.18}
\end{equation}
where $s$ can be chosen such that the sufficient condition have to
hold. Note that from the expansions given in (\ref{eq:beta_k exp}-\ref{eq:3c})
and from (\ref{eq:3.18}) we have
\[
\omega_{k}^{NS7}=\begin{cases}
d_{k}+O\left(\Delta x^{6}\right)^{s}, & \text{ if \,}f_{j+\frac{1}{2}}^{(1)}\neq0,\\
d_{k}+O\left(\Delta x^{4}\right)^{s}, & \text{ if \,}f_{j+\frac{1}{2}}^{(1)}=0,f_{j+\frac{1}{2}}^{(2)}\neq0,\\
d_{k}+O\left(\Delta x^{2}\right)^{s}, & \text{ if }\,f_{j+\frac{1}{2}}^{(1)}=0,f_{j+\frac{1}{2}}^{(2)}=0,f_{j+\frac{1}{2}}^{(3)}\neq0.
\end{cases}
\]
Clearly, the sufficient condition (\ref{eq:17}) is satisfies for
the WENO-NS7 scheme under the following conditions:
\begin{enumerate}
\item $s\geq1$ if $f_{j+\frac{1}{2}}^{(1)}\neq0,$ or if $f_{j+\frac{1}{2}}^{(1)}=0\text{ and }f_{j+\frac{1}{2}}^{(2)}\neq0.$
\item $s\geq2$ if $f_{j+\frac{1}{2}}^{(1)}=0,\;f_{j+\frac{1}{2}}^{(2)}=0,\;f_{j+\frac{1}{2}}^{(3)}\neq0$.
\end{enumerate}
For numerical verification, $s$ value is taken as $2.$ With these
developments, in the next section we'll test the WENO-NS7 scheme for
various examples.

\section{Numerical results}

\hspace{0.6 cm}Let's denote the system (\ref{eq:1}) by
\[
\frac{du}{dt}=L(u),
\]
where $L(u)$ is an approximation to the derivative $-f(u)_{x}.$
In section 2, we have obtained higher order reconstruction for the
flux function which is defined in (\ref{eq:2}). To evolve the solution
in time, strong-stability-preserving Runge\textendash Kutta $(SSPRK)$
algorithm is used, whose detailed description is in article \cite{Gottilieb}.
The choice of this time integration is to ensure that the order of
accuracy for the time evaluation matches with that of the spatial
order of accuracy. For linear problems, $m-$stage linear $SSPRK$
method, which is of $(m-1)^{th}$ order is used in
the following examples. The $lSSPRK(m,m-1)$ method
is
\begin{align*}
 & u^{(0)}=u^{n},\\
 & u^{(i)}=u^{(i-1)}+\frac{1}{2}\Delta tL\left(u^{(i-1)}\right),i=1,...,m-1\\
 & u^{(m)}=u^{n+1}=\sum_{k=0}^{m-2}\alpha_{m,k}u^{(k)}+\alpha_{m,m-1}\left(u^{(m-1)}+\frac{1}{2}\Delta tL\left(u^{(m-1)}\right)\right).
\end{align*}
For the seventh order the coefficients $\alpha_{m,k},(m=8)$ are given
as
\[
\alpha_{8,0}=\frac{2}{15},\,\alpha_{8,1}=\frac{2}{7},\,\alpha_{8,2}=\frac{2}{9},\,\alpha_{8,3}=\frac{4}{15},\,\alpha_{8,4}=0,\,\alpha_{8,5}=\frac{4}{45},\,\alpha_{8,6}=0,\,\alpha_{8,7}=\frac{1}{315}.
\]
This method will not attain $O\left(\Delta t{}^{(m)}\right)$ for
nonlinear problems. So, for nonlinear problems a fourth order nonlinear
version of $SSPRK(5,4)$ is used with the stability condition $CFL\leq1,$
where $CFL=\underset{j}{\max}\left\{ S_{j}^{n}\frac{\Delta t}{\Delta x}\right\} .$
Here $S_{j}^{n}$ is the maximum propagation speed in $I_{j}$ at
time level $t^{n}$. The fourth order $SSPRK(5,4)$ method is given
as
\begin{align*}
u^{(1)} & =u^{n}+0.39175222\Delta tL\left(u^{n}\right),\\
u^{(2)} & =0.444370494\Delta u^{n}+0.5556295u^{(1)}+0.36841059\Delta tL\left(u^{(1)}\right),\\
u^{(3)} & =0.62010185\Delta u^{n}+0.379898u^{(2)}+0.25189177\Delta tL\left(u^{(2)}\right),\\
u^{(4)} & =0.17807995u^{n}+0.82192004u^{(3)}+0.5449747\Delta tL\left(u^{(3)}\right),\\
u^{n+1} & =0.517231u^{(2)}+0.0960597u^{(3)}+0.3867086u^{(4)}+0.063692\Delta tL\left(u^{(3)}\right)+0.22600748\Delta tL\left(u^{(4)}\right).
\end{align*}
For numerical comparison of the proposed WENO-NS7 scheme, the numerical
results are presented along with the numerical results of seventh
order WENO-BS \cite{balsara and shu8} and seventh order WENO-Z \cite{W S DON }
schemes in the following sections. These results are mostly comparable with WENO-BS
scheme in compare to other seventh-order WENO schemes. For completeness these schemes
are briefly described in the appendix.

\subsection{Scalar test problems}

\hspace{0.6 cm}To verify the numerical accuracy and convergence of
the scheme examples pertaining to transport and Burger's equations
with various initial profiles are considered. Some of these initial
profiles contain jump discontinuity and in some cases, the solution
in time leads to shocks. Lax-Friedrich's flux splitting technique
is used in (\ref{eq:3-1}). For the scheme WENO-BS, $\epsilon=10^{-6}$
and for WENO-Z and WENO-NS7 schemes $\epsilon=10^{-40}$ is taken
along with the CFL number $0.5$. The parameters in (\ref{eq:beta_k})
are fixed as $\xi_{1}=0.1$ and $\xi_{2}=1$ for linear test cases.

\subsubsection{Behaviors of new weights}
\begin{figure}[H]
\includegraphics[width=16cm,height=12cm]{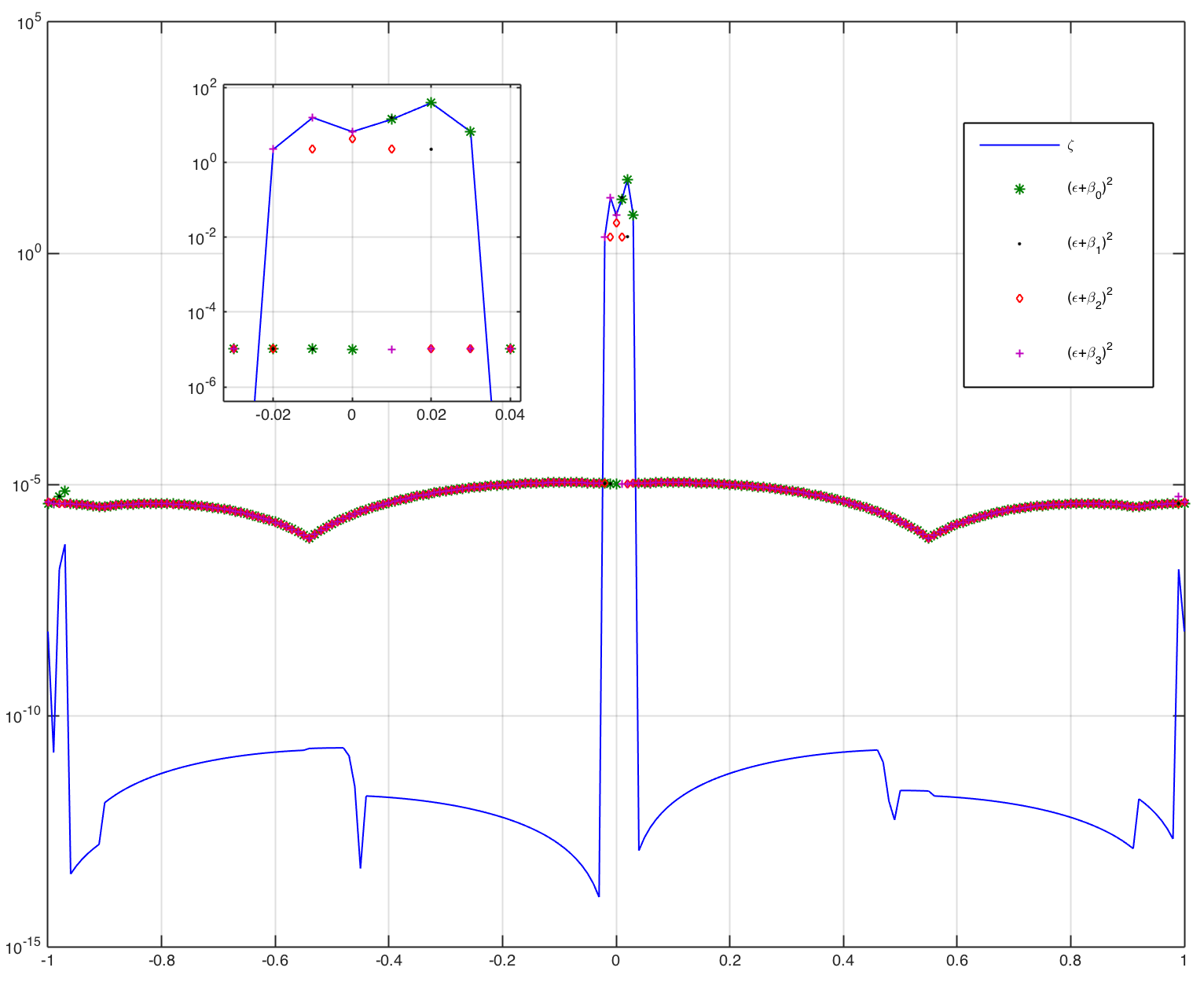}
\caption{\label{fig:SI}The values of the smoothness indicators $\beta_{k}$
to initial data ($\ref{eq:4.8})$ by WENO-NS7 scheme}
\end{figure}

\begin{figure}[H]
\centering{}$\qquad$\includegraphics[width=14cm,height=10cm]{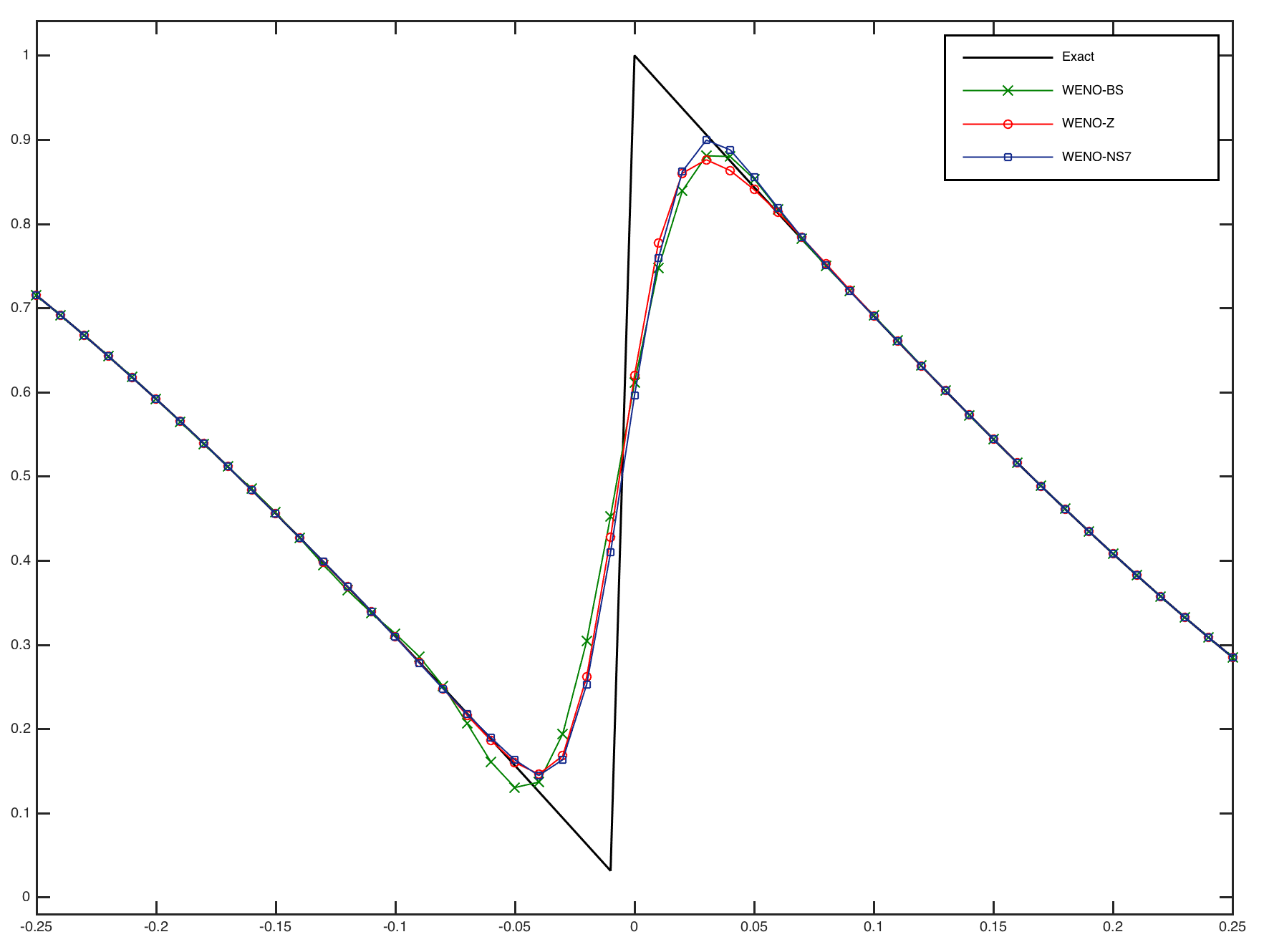}\caption{\label{fig:NS1}Numerical solution of linear advection equation with
the initial condition ($\ref{eq:4.8}$)}
\end{figure}

\begin{figure}[H]
\centering{}\includegraphics[width=20cm,height=7cm]{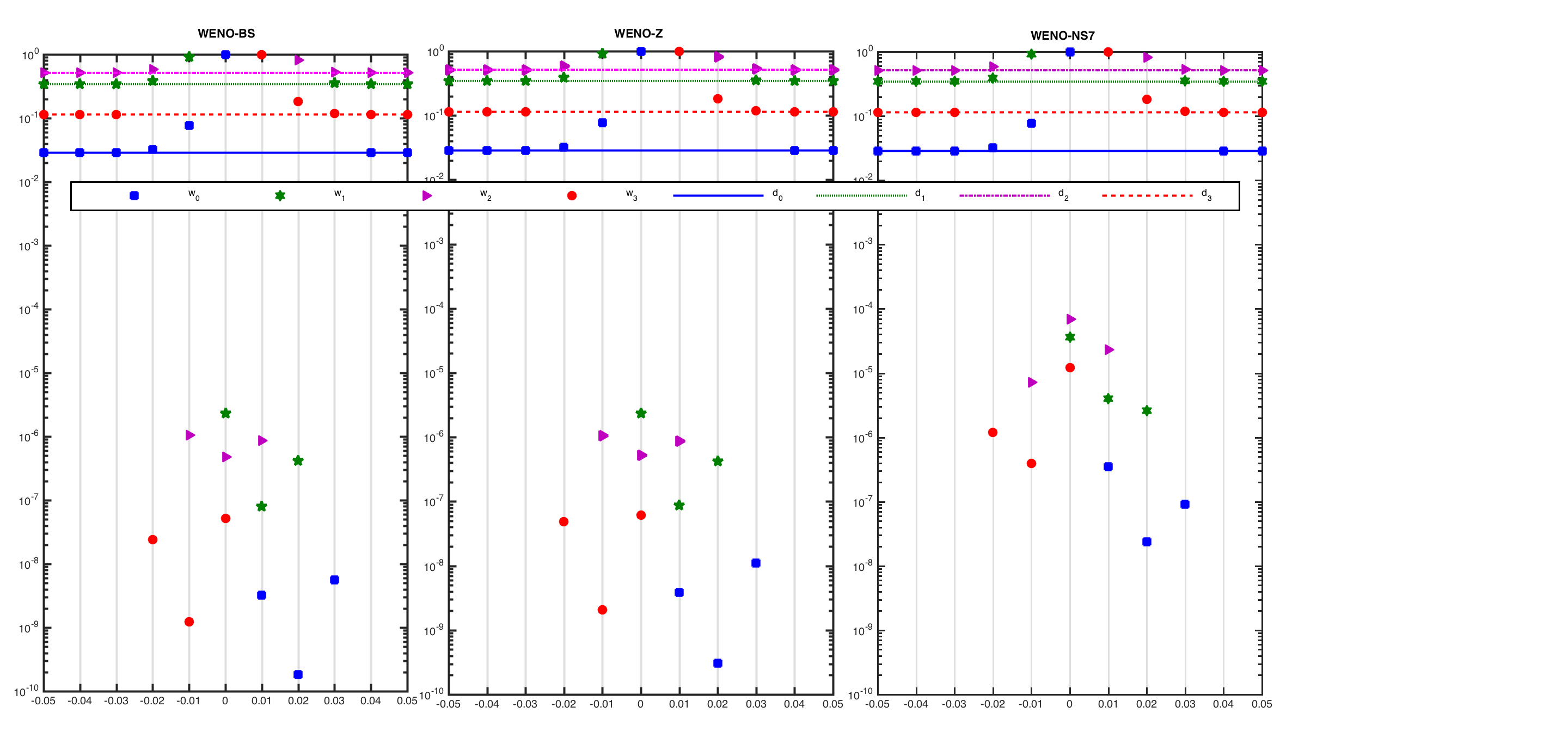}\caption{\label{fig:Weights}The distribution of ideal weights $d_{k}$ and
the weights $\omega_{k},$$k=0,1,2,3$. }
\end{figure}

\textbf{Example 1:}
For linear advection equation
\[
u_{t}+u_{x}=0,-1\leq x\leq1,t>0,
\]
let the initial condition be
\begin{equation}
u(x,0)=u_{0}(x)=\begin{cases}
-\text{sin}(\pi x)-\frac{1}{2}x^{3} & \text{ for }\,\,-1\leq x<0,\\
-\text{sin}(\pi x)-\frac{1}{2}x^{3}+1 & \text{ for }\,\,0\leq x<1,
\end{cases}\label{eq:4.8}
\end{equation}
which is a piecewise continuous function with jump discontinuity at
$x=0.$
The behavior of the smoothness indicators $\beta_{k}\,(k=0,1,2,3)$
and the global smoothness indicator $\zeta$ for initial time, $t=0$
is displayed in figure \ref{fig:SI} for the proposed scheme WENO-NS7.
The approximate solution is computed with uniform discretization of
the spatial domain with the step size, $\triangle x=0.01$ up to time
$t=8$ with the scheme WENO-NS7 along with WENO-BS and WENO-Z schemes,
these solutions are plotted in figure \ref{fig:NS1} against the exact
solution. It can be observed from the plot that the proposed scheme
performs better than other schemes near the jump discontinuity. Figure
\ref{fig:Weights} displays the behavior of the weights $\omega_{k}$
along with the ideal weights $d_{k}.$

\subsubsection{Accuracy, at critical points}
\textbf{Example 2:}
Consider the linear transport equation
\begin{equation}
u_{t}+u_{x}=0,\,\,-1\leq x\leq1,\,\,t>0,\label{eq:4.4}
\end{equation}
with the periodic boundary conditions up to time $t=2.$ Three different
initial conditions are considered, each of them is a special case
to test the convergence analysis.\\
$\textbf{Case 1:}$ The smooth initial condition
\begin{equation}
u(x,0)=\text{sin}(\pi x),\label{eq:4.5}
\end{equation}
is taken in this case to verify the order of convergence.

In table \ref{table:1}, the $L_{1}$ and $L_{\infty}$ errors along
with the numerical order of accuracy are given for WENO-BS, WENO-Z
and WENO-NS7 schemes. It has been observed that the proposed scheme
attains the desired order of accuracy and very much efficient than
WENO-BS scheme but has the same accuracy as WENO-Z scheme.
\begin{table}[h]
{\scriptsize{}$\qquad\qquad$}%
\begin{tabular}{rcccccc}
\hline
 & {WENO-BS} &  & WENO-Z &  & WENO-NS7 & \tabularnewline
\hline
N & $L_{1}\text{ error}$ & $L_{1}\text{ order}$ & $L_{1}\text{ error}$ & $L_{1}\text{ order}$ & $L_{1}\text{ error}$ & $L_{1}\text{ order}$\tabularnewline
\hline
10 & 3.9281e-03 & \textemdash{} & 6.5177e-04 & \textemdash{} & 9.9734e-04 & \textemdash{}\tabularnewline
20 & 8.5335e-05 & 5.52 & 4.2296e-06 & 7.26 & 5.6871e-06 & 7.45\tabularnewline
40 & 1.2222e-06 & 6.12 & 3.3547e-08 & 6.98 & 3.3551e-08 & 7.40\tabularnewline
80 & 1.8418e-08 & 6.05 & 2.6304e-10 & 6.99 & 2.6304e-10 & 6.99\tabularnewline
160 & 2.7931e-10 & 6.04 & 2.0638e-12 & 6.99 & 2.0637e-12 & 6.99\tabularnewline
\hline
 & $L_{\infty}\text{ error}$ & $L_{\infty}\text{ order}$ & $L_{\infty}\text{ error}$ & $L_{\infty}\text{ order}$ & $L_{\infty}\text{ error}$ & $L_{\infty}\text{ order}$\tabularnewline
\hline
10 & 8.5830e-03 & \textemdash{} & 1.4372e-03 & \textemdash{} & 2.0586e-03 & \textemdash{}\tabularnewline
20 & 2.2095e-04 & 5.28 & 6.7946e-06 & 7.72 & 1.5028e-05 & 7.09\tabularnewline
40 & 5.3163e-06 & 5.38 & 5.2546e-08 & 7.01 & 6.8261e-08 & 7.78\tabularnewline
80 & 1.4576e-07 & 5.18 & 4.1283e-10 & 6.99 & 4.3546e-10 & 7.29\tabularnewline
160 & 4.0863e-09 & 5.15 & 3.2415e-12 & 6.99 & 3.2736e-12 & 7.05\tabularnewline
\hline
\end{tabular}

\caption{\label{table:1}$L_{1}\text{and }L_{\infty}-$ errors of example (\ref{eq:4.4})
with initial condition (\ref{eq:4.5})}
\end{table}
\\
$\textbf{Case 2:}$ In this case, the initial condition is chosen
as
\begin{equation}
u(x,0)=\text{sin}(\pi x-\frac{1}{\pi}\text{sin}(\pi x)),\label{eq:4.6}
\end{equation}
which contains first-order critical point, that is, $u_{x}=0$ but
note that $u_{xxx}\neq0.$

The $L_{1}$ and $L_{\infty}$ errors along with the numerical order
of accuracy are provided in table \ref{table:2} for WENO-BS, WENO-Z
and WENO-NS7 schemes. The proposed scheme WENO-NS7 achieves the desired
order of accuracy.
\begin{table}[h]
{\scriptsize{}$\qquad\qquad$}%
\begin{tabular}{rcccccc}
\hline
 & {WENO-BS} &  & WENO-Z &  & WENO-NS7 & \tabularnewline
\hline
N & $L_{1}\text{ error}$ & $L_{1}\text{ order}$ & $L_{1}\text{ error}$ & $L_{1}\text{ order}$ & $L_{1}\text{ error}$ & $L_{1}\text{ order}$\tabularnewline
\hline
10 & 2.9823e-02 & \textemdash \textemdash{} & 3.0006e-02 & \textemdash \textendash{} & 1.3986e-02 & \textemdash \textendash{}\tabularnewline
20 & 7.5922e-04 & 5.29 & 5.0945e-04 & 5.88 & 3.0947e-04 & 5.49\tabularnewline
40 & 8.3317e-06 & 6.51 & 2.6348e-06 & 7.59 & 2.6567e-06 & 6.86\tabularnewline
80 & 6.9845e-08 & 6.89 & 2.1488e-08 & 6.93 & 2.1524e-08 & 6.94\tabularnewline
160 & 7.1684e-10 & 6.60 & 1.6933e-10 & 6.99 & 1.6934e-10 & 6.99\tabularnewline
\hline
 & $L_{\infty}\text{ error}$ & $L_{\infty}\text{ order}$ & $L_{\infty}\text{ error}$ & $L_{\infty}\text{ order}$ & $L_{\infty}\text{ error}$ & $L_{\infty}\text{ order}$\tabularnewline
\hline
10 & 7.1325e-02 & \textemdash \textemdash{} & 7.0754e-02 & \textemdash \textemdash{} & 3.6883e-02 & \textemdash \textemdash{}\tabularnewline
20 & 2.2429e-03 & 4.99 & 1.3617e-03 & 5.70 & 9.2089e-04 & 5.32\tabularnewline
40 & 3.6731e-05 & 5.93 & 7.9578e-06 & 7.41 & 7.9697e-06 & 6.85\tabularnewline
80 & 5.8755e-07 & 5.96 & 6.6439e-08 & 6.90 & 6.6272e-08 & 6.91\tabularnewline
160 & 1.1019e-08 & 5.73 & 5.2713e-10 & 6.97 & 5.2711e-10 & 6.97\tabularnewline
\hline
\end{tabular}

\caption{\label{table:2}$L_{1}\text{and }L_{\infty}-$ errors of example (\ref{eq:4.4})
with initial condition (\ref{eq:4.6})}
\end{table}
\\
$\textbf{Case 3:}$ The initial condition
\begin{equation}
u(x,0)=\text{sin}(\pi x)^{3},\label{eq:4.7}
\end{equation}
has the nature, $u_{x}=0,\:u_{xx}=0$ but $u_{xxx}\neq0.$

In table \ref{table:3}, the $L_{1}$ and $L_{\infty}$ errors are
tabulated along with the numerical order of accuracy for the schemes
WENO-BS, WENO-Z and WENO-NS7. In this case too, the proposed scheme
has the desired order of convergence.
\begin{table}[H]
{\scriptsize{}$\qquad\qquad$}%
\begin{tabular}{rcccccc}
\hline
 &{WENO-BS} &  & WENO-Z &  & WENO-NS7 & \tabularnewline
\hline
N & $L_{1}\text{ error}$ & $L_{1}\text{ order}$ & $L_{1}\text{ error}$ & $L_{1}\text{ order}$ & $L_{1}\text{ error}$ & $L_{1}\text{ order}$\tabularnewline
\hline
10 & 2.1618e-01 & \textemdash \textemdash{} & 2.1128e-01 & \textemdash \textendash{} & 1.9695e-01 & \textemdash \textendash{}\tabularnewline
20 & 1.8858e-02 & 3.52 & 8.4032e-03 & 4.65 & 8.6604e-03 & 4.50\tabularnewline
40 & 4.1300e-04 & 5.51 & 6.4354e-05 & 7.02 & 8.6452e-05 & 6.64\tabularnewline
80 & 5.9000e-06 & 6.13 & 4.2948e-07 & 7.22 & 4.1937e-07 & 7.68\tabularnewline
160 & 6.6306e-08 & 6.47 & 3.3639e-09 & 6.99 & 3.3582e-09 & 6.96\tabularnewline
\hline
 & $L_{\infty}\text{ error}$ & $L_{\infty}\text{ order}$ & $L_{\infty}\text{ error}$ & $L_{\infty}\text{ order}$ & $L_{\infty}\text{ error}$ & $L_{\infty}\text{ order}$\tabularnewline
\hline
10 & 2.9623e-01 & \textemdash \textemdash{} & 2.9814e-01 & \textemdash \textemdash - & 3.2319e-01 & \tabularnewline
20 & 3.8372e-02 & 2.95 & 1.4710e-02 & 4.34 & 1.7309e-02 & 4.22\tabularnewline
40 & 1.0319e-03 & 5.21 & 1.4691e-04 & 6.64 & 2.1456e-04 & 6.33\tabularnewline
80 & 2.3364e-05 & 5.46 & 6.6758e-07 & 7.78 & 1.0238e-06 & 7.71\tabularnewline
160 & 6.1191e-07 & 5.25 & 5.2805e-09 & 6.98 & 5.5897e-09 & 7.51\tabularnewline
\hline
\end{tabular}

\caption{\label{table:3}$L_{1}\text{and }L_{\infty}-$ errors of example (\ref{eq:4.4})
with initial condition (\ref{eq:4.7})}
\end{table}

\begin{figure}[H]
\includegraphics[width=18cm,height=8cm]{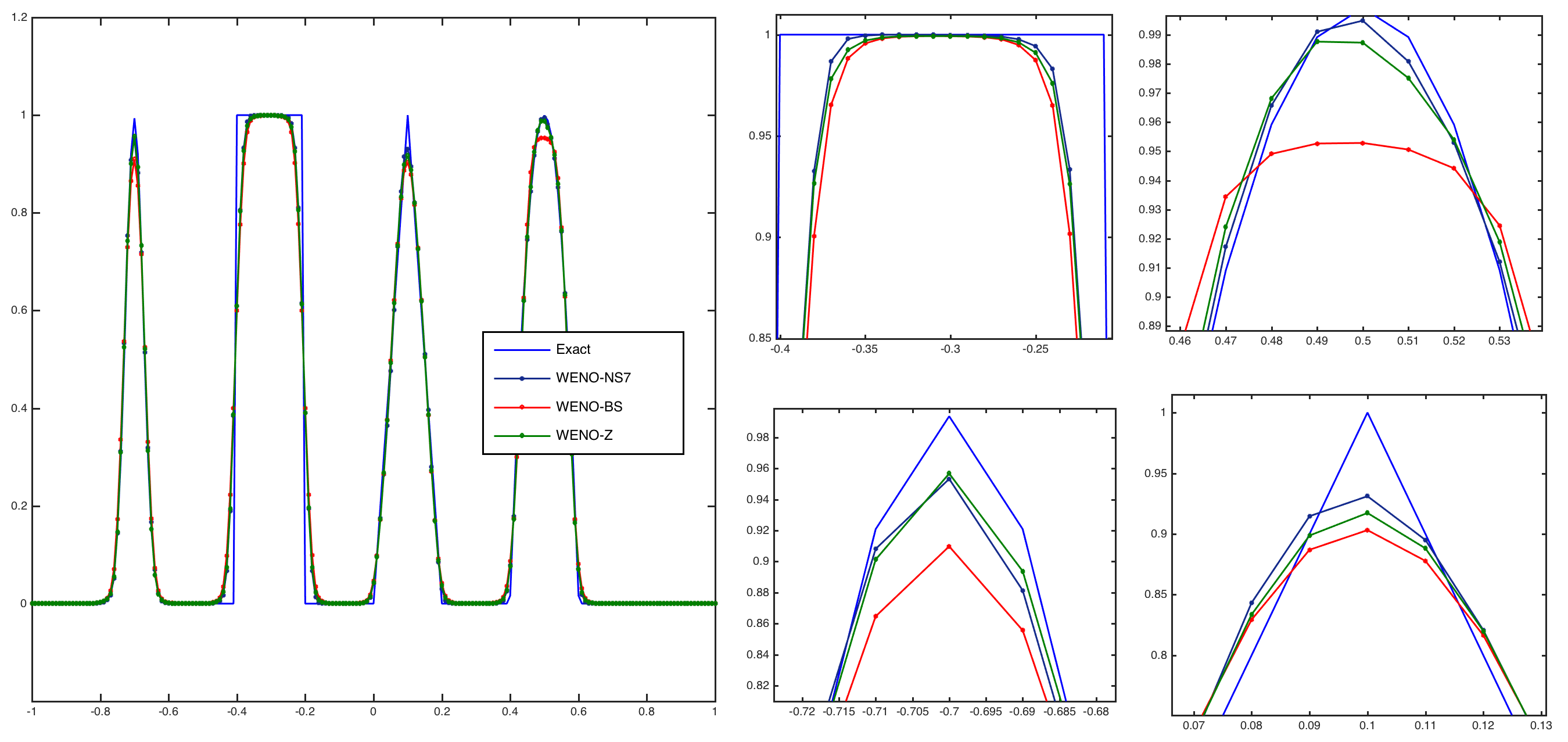}\caption{\label{fig:NS2}Numerical solution of linear advection equation with
the initial condition ($\ref{eq:4.9-1}$)}
\end{figure}

\textbf{Example 3:}
Consider the linear advection equation $u_{t}+u{}_{x}=0$ on an interval
$[-1,1]$ with the initial condition
\begin{equation}
u(x,0)=\begin{cases}
\frac{1}{6}[G(x,z-\delta)+G(x,z+\delta)+4G(x,z)], & -0.8\leq x\leq-0.6\\
1, & -0.4\leq x\leq-0.2\\
1-|10(x-0.1)|, & 0\leq x\leq0.2\\
\frac{1}{6}[F(x,a-\delta)+F(x,a+\delta)+4F(x,a)], & 0.4\leq x\leq0.6\\
0, & \mathrm{otherwise}
\end{cases}\label{eq:4.9-1}
\end{equation}
where $G(x,z)=\textrm{exp}\left(-\beta(x-z)^{2}\right)$, $F(x,a)=\left\{ \textrm{max}\left(1-\alpha^{2}(x-a)^{2},0\right)\right\} {}^{\frac{1}{2}}$,
$a=0.5$, $z=-0.7$, $\delta=0.005$, $\alpha=10$ and  $\beta=\left(\frac{\log(2)}{36\delta^{2}}\right)$.
This initial condition consists of several shapes that contains the
combination of Gaussian, a square wave, a sharp triangle wave and
an half ellipse, which are difficult to resolve by the numerical schemes,
so it turns out to be a typical test case for numerical verification.
The solution is computed with periodic boundary conditions on the
mesh of 200 cells up to time $t=8$ by the schemes WENO-NS7, WENO-BS
and WENO-Z. These solutions are plotted in figure \ref{fig:NS2} against
the exact solution, from the figure one can observe that WENO-NS7
scheme performs better than WENO-BS and WENO-Z schemes.
\subsubsection{Test with shocks}
\textbf{Example 4:}
Consider the inviscid Burger's equation
\begin{equation}
u_{t}+\left(\frac{u^{2}}{2}\right){}_{x}=0,\,\,-1\leq x\leq1,\,\,t>0,\label{eq:4.9}
\end{equation}
 with the following initial conditions.

$\textbf{Case 1:}$ The initial condition
\begin{equation}
u_{0}(x)=-\text{sin}(\pi x),\label{eq:4.11-1}
\end{equation}
produces a steady shock at the point $x=0$ as the time progresses
from $t=0$ to $t=1.5.$

\begin{figure}[H]
\centering{}\includegraphics[width=14cm,height=7cm]{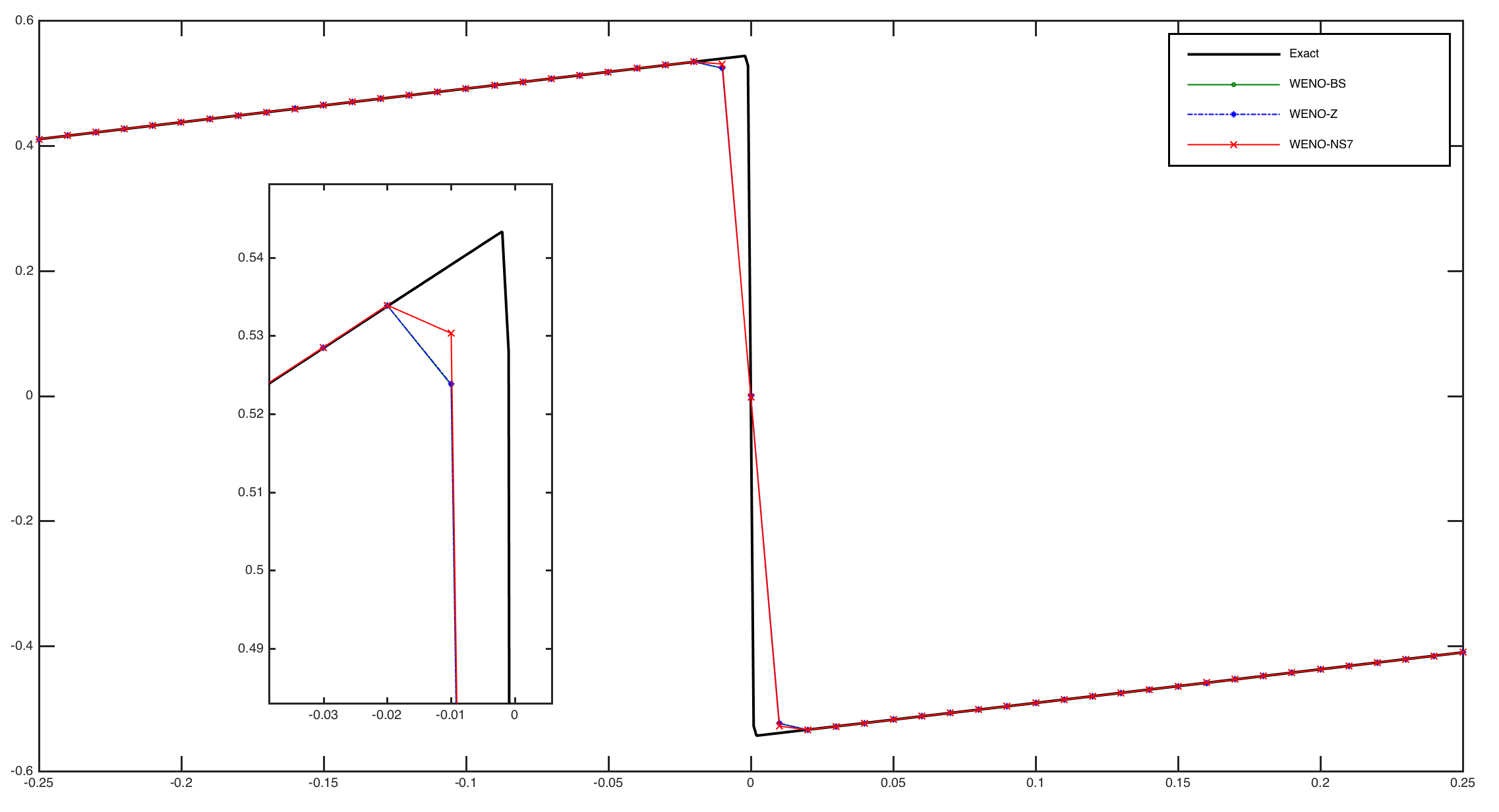}
\caption{\label{fig:Burg1}Approximate solution of (\ref{eq:4.9}) with initial
condition (\ref{eq:4.11-1}) and its zoomed region. }
\end{figure}

The numerical solution is computed with the periodic boundary conditions,
on the mesh of $200$ grid cells. The approximate solutions are plotted
in figure \ref{fig:Burg1} against the reference (exact) solution.
It is noted from this figure that the WENO-NS7 scheme captures the
shock very well in compare WENO-BS and WENO-Z schemes.

\begin{figure}[H]
\centering{}\includegraphics[width=14cm,height=7cm]{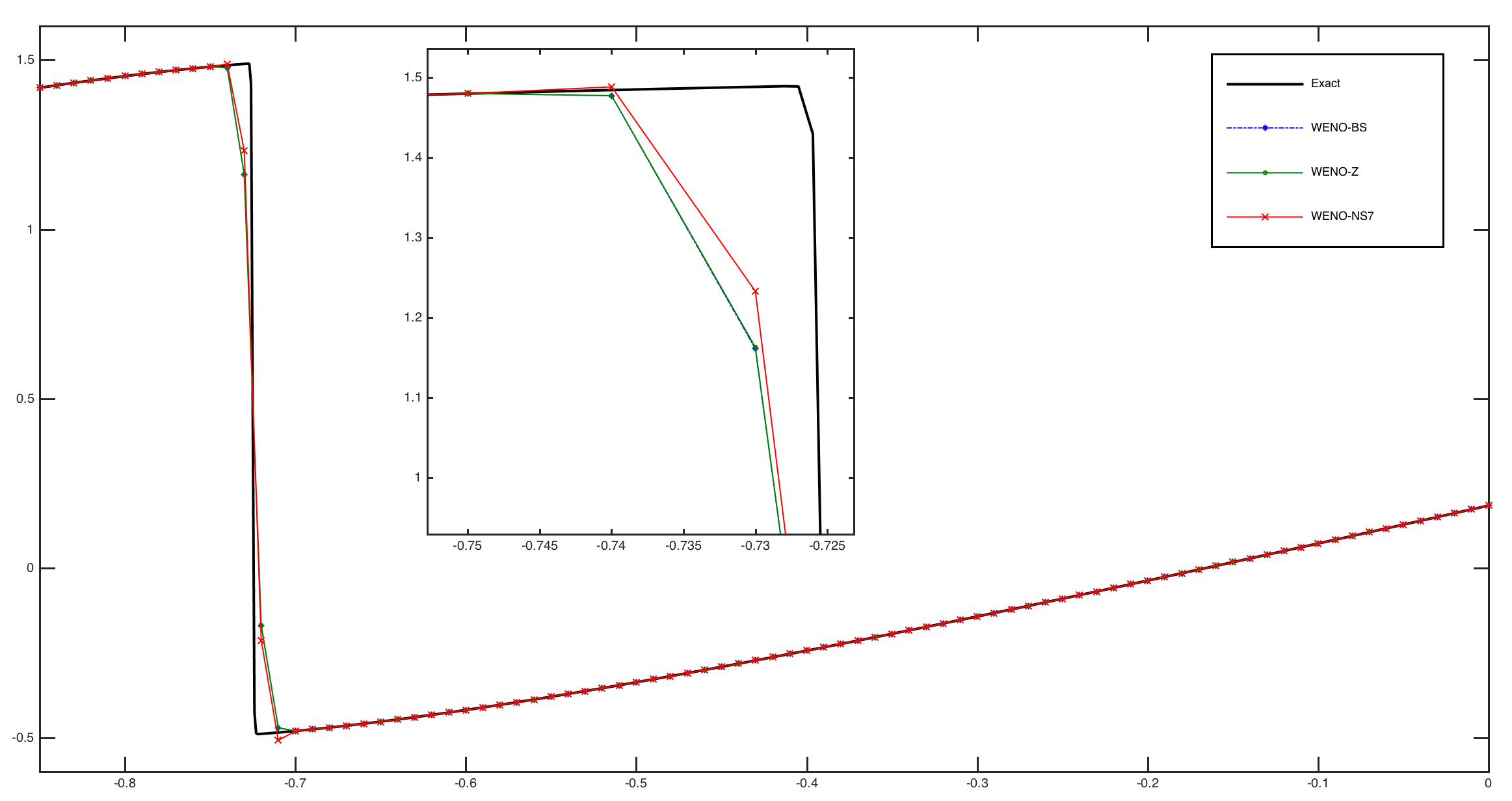}\caption{\label{fig:Burg2}Approximate solution of (\ref{eq:4.9}) with initial
condition (\ref{eq:4.12}) and its zoomed region.}
\end{figure}

$\textbf{Case 2:}$ Consider the initial condition
\begin{equation}
u_{0}(x)=\frac{1}{2}+\text{sin}(\pi x),\label{eq:4.12}
\end{equation}
for the equation (\ref{eq:4.9}), which generates a moving shock-wave
at time $t=0.55.$ The numerical results for the WENO-NS7 along with
the WENO-BS and WENO-Z schemes are plotted in the figure \ref{fig:Burg2}
on the computational grid of $200$ cells against the reference solution.
Here too, WENO-NS7 scheme captures the shock better than other schemes.
\textbf{Remark}
We've chosen the parameters $\xi_{1}=0.1$ and $\xi_{2}=0.3$ to evaluate
the equation (\ref{eq:beta_k}) for the solution of (\ref{eq:4.9}).

\subsection{Euler equations in one space dimension}

\hspace{0.6 cm}The one-dimensional system of Euler equations are
given by
\begin{equation}
\begin{pmatrix}\rho\\
\rho u\\
E
\end{pmatrix}_{t}+\begin{pmatrix}\rho u\\
\rho u^{2}+p\\
u(E+p)
\end{pmatrix}_{x}=0\label{eq:60}
\end{equation}
where $\rho,u,E,p$ are the density, velocity, total energy and pressure
respectively.

The system (\ref{eq:60}) represents the conservation of mass, momentum
and total energy. The total energy for an ideal polytropic gas is
defined as
\begin{equation}
E=\frac{p}{\gamma-1}+\frac{1}{2}\rho u^{2},
\end{equation}
where $\gamma$ is the ratio of specific heats and its value is taken
as $\gamma=1.4.$

For the numerical evaluation of this system by WENO-NS7 scheme, in
the following examples the parameters in (\ref{eq:3-1}) are taken
as $\xi_{1}=0.3$ and $\xi_{2}=0.3.$

\subsubsection{Shock tube problems}

\begin{figure}[H]
\centering{}\includegraphics[width=14cm,height=6cm]{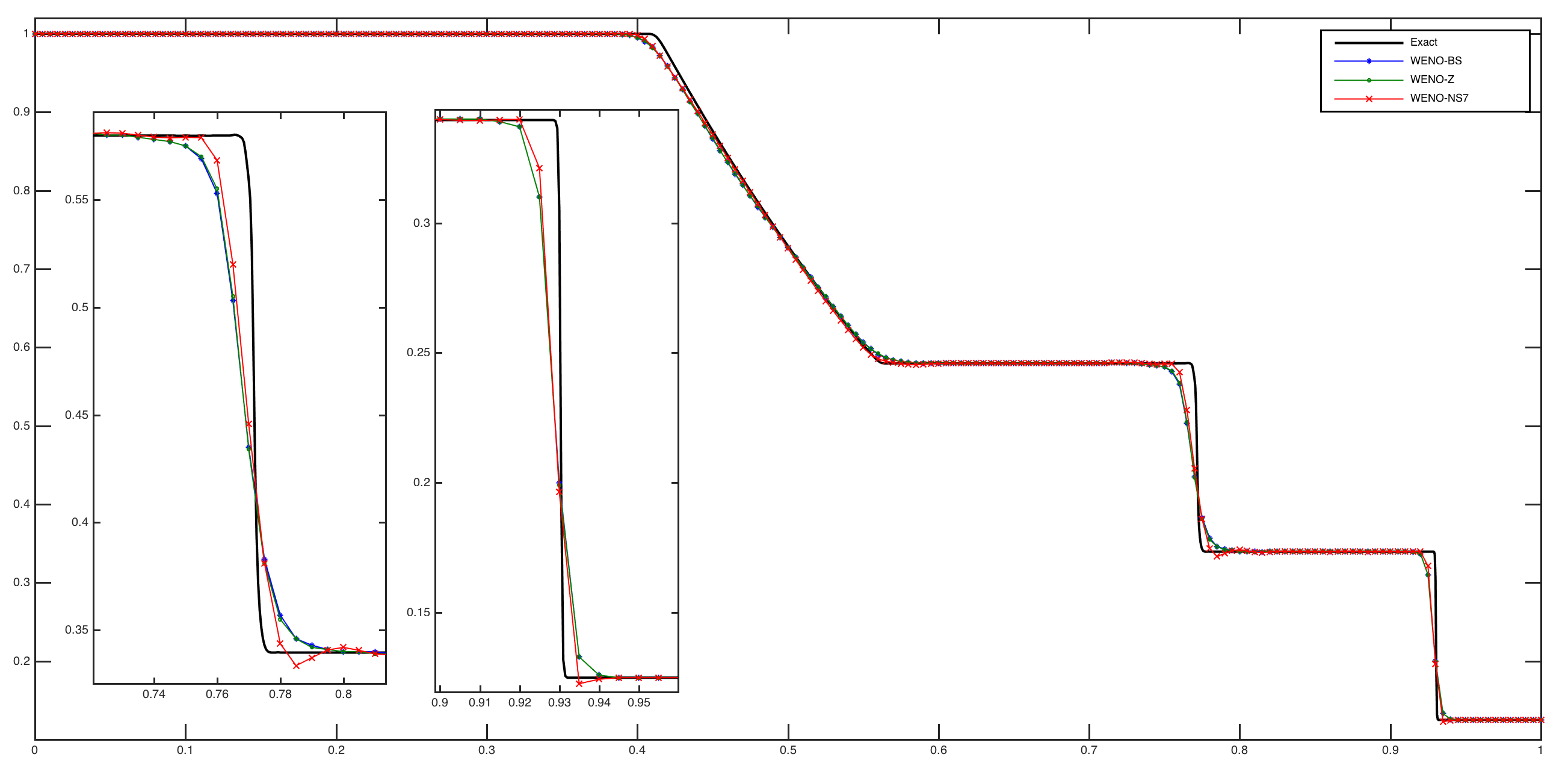}\caption{\label{fig:Den-sod1}Density profile for Sod's shock tube problem.}
\end{figure}

\begin{figure}[H]
\centering{}\includegraphics[width=14cm,height=6cm]{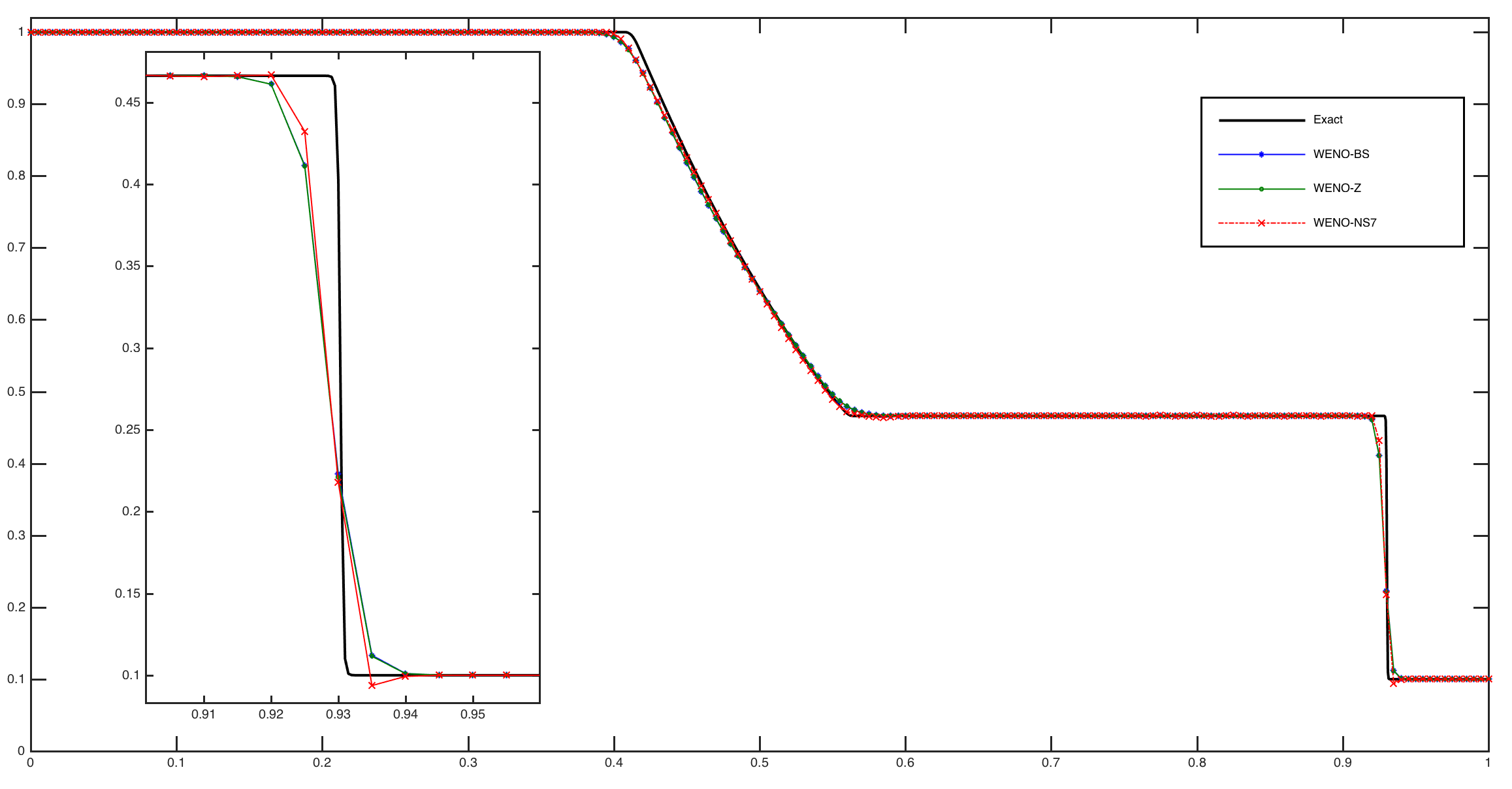}\caption{\label{fig:Pres-sod1}Pressure profile of Sod's shock tube problem. }
\end{figure}

\begin{figure}[H]
\centering{}\includegraphics[width=14cm,height=6cm]{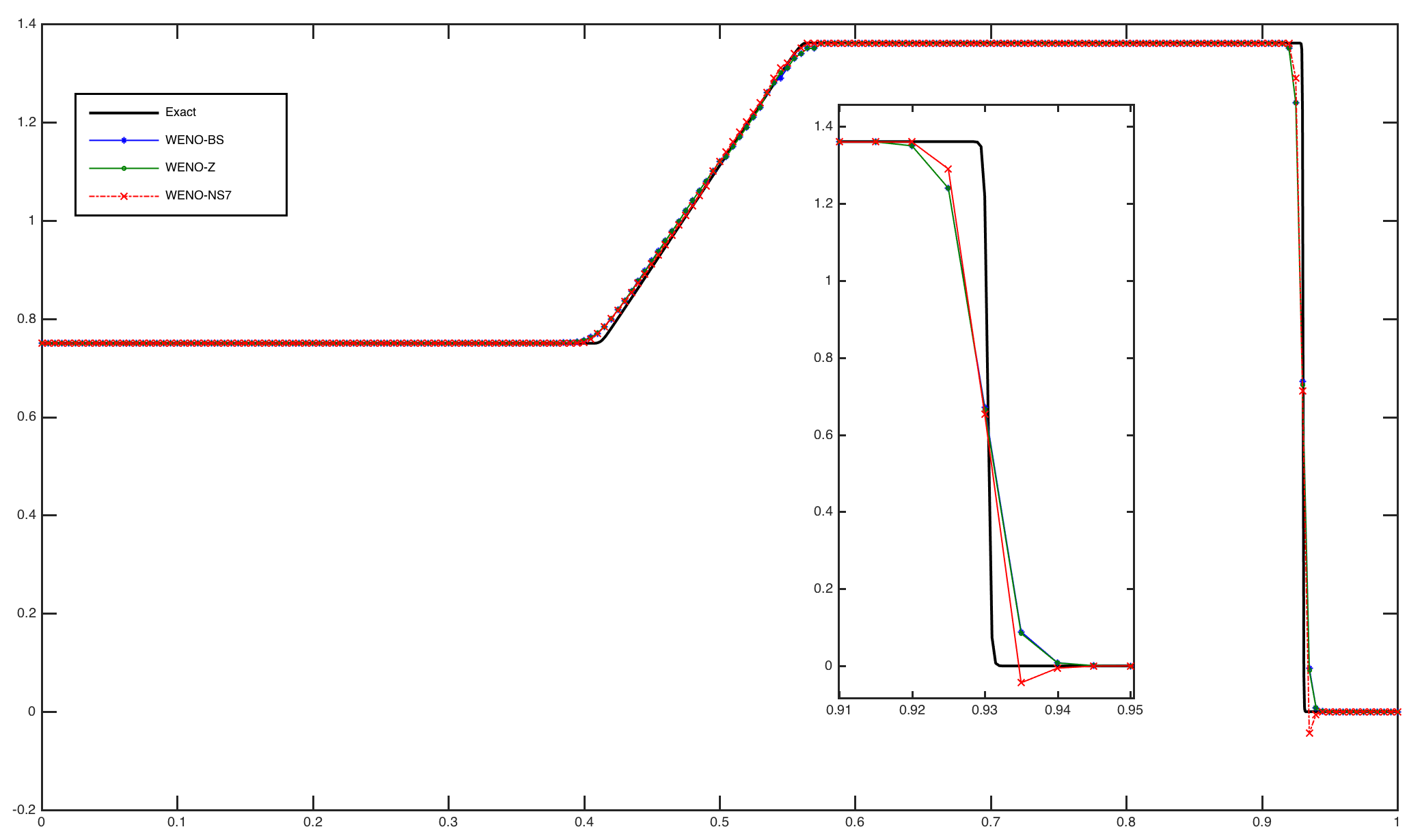}\caption{\label{fig:velo-sod1}Velocity profile of Sod's shock tube problem. }
\end{figure}

Consider the system (\ref{eq:60}) with the following well-known initial
conditions:\\
\textbf{Example 5:}
The Riemann data
\begin{equation}
(\rho,u,p)=\begin{cases}
(1.0,0.75,1.0), & \text{ if }0\leq x\leq0.5,\\
(0.125,0.0,0.1), & \text{ if }0.5\leq x\leq1,
\end{cases}\label{58}
\end{equation}
is a modified version of Sod's shock tube problem \cite{Sod}, its
solution contains a right shock wave, a right traveling contact wave
and a left sonic rarefaction wave.
The numerical results are obtained on a grid of $200$ cells for the
spatial domain $0\leq x\leq1$ up to $t=0.2$ with the transmissive
boundary conditions. The density, pressure and velocity profiles are
plotted in figures \ref{fig:Den-sod1}, \ref{fig:Pres-sod1} and \ref{fig:velo-sod1}
respectively, against the reference solution and the diagrams contain
the enlarged portions of shock and contact discontinuity regions.
From these figures, one can conclude that WENO-NS7 scheme performs
better in compare to WENO-BS and WENO-Z schemes.
\textbf{Example 6:}
The initial condition
\begin{equation}
(\rho,u,p)=\begin{cases}
(0.445,0.698,3.528), & \text{ if }-5\leq x<0,\\
(0.500,0.000,0.571), & \text{ if }0\leq x\leq5,
\end{cases}\label{59}
\end{equation}
is due to Lax \cite{PDLax}. This problem reveals the shock capturing
capabilities of the numerical scheme.

\begin{figure}[H]
\centering{}\includegraphics[width=14cm,height=7cm]{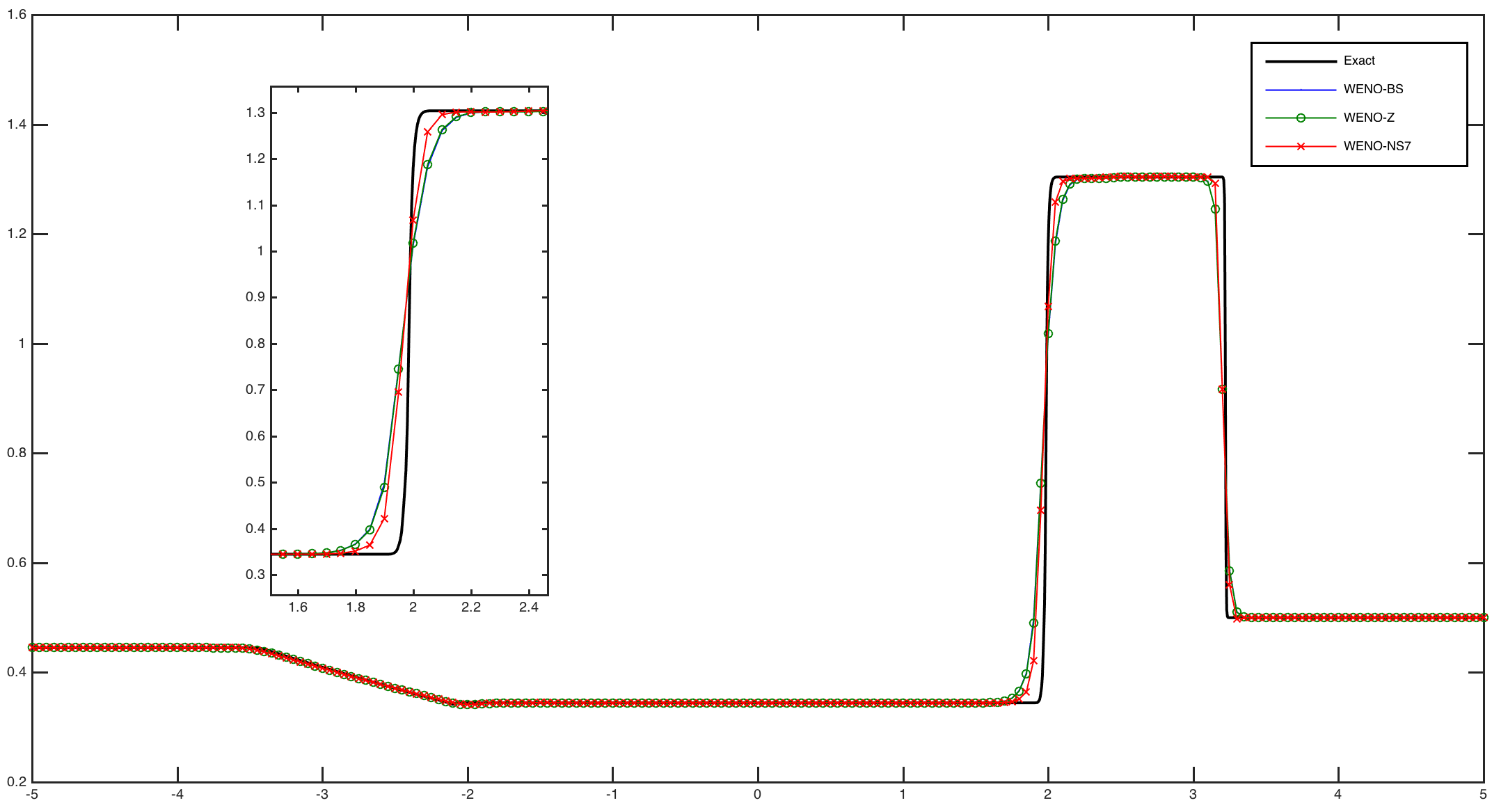}\caption{\label{fig:Den-Lax}Density profile for Lax initial condition.}
\end{figure}

\begin{figure}[H]
\centering{}\includegraphics[width=14cm,height=6cm]{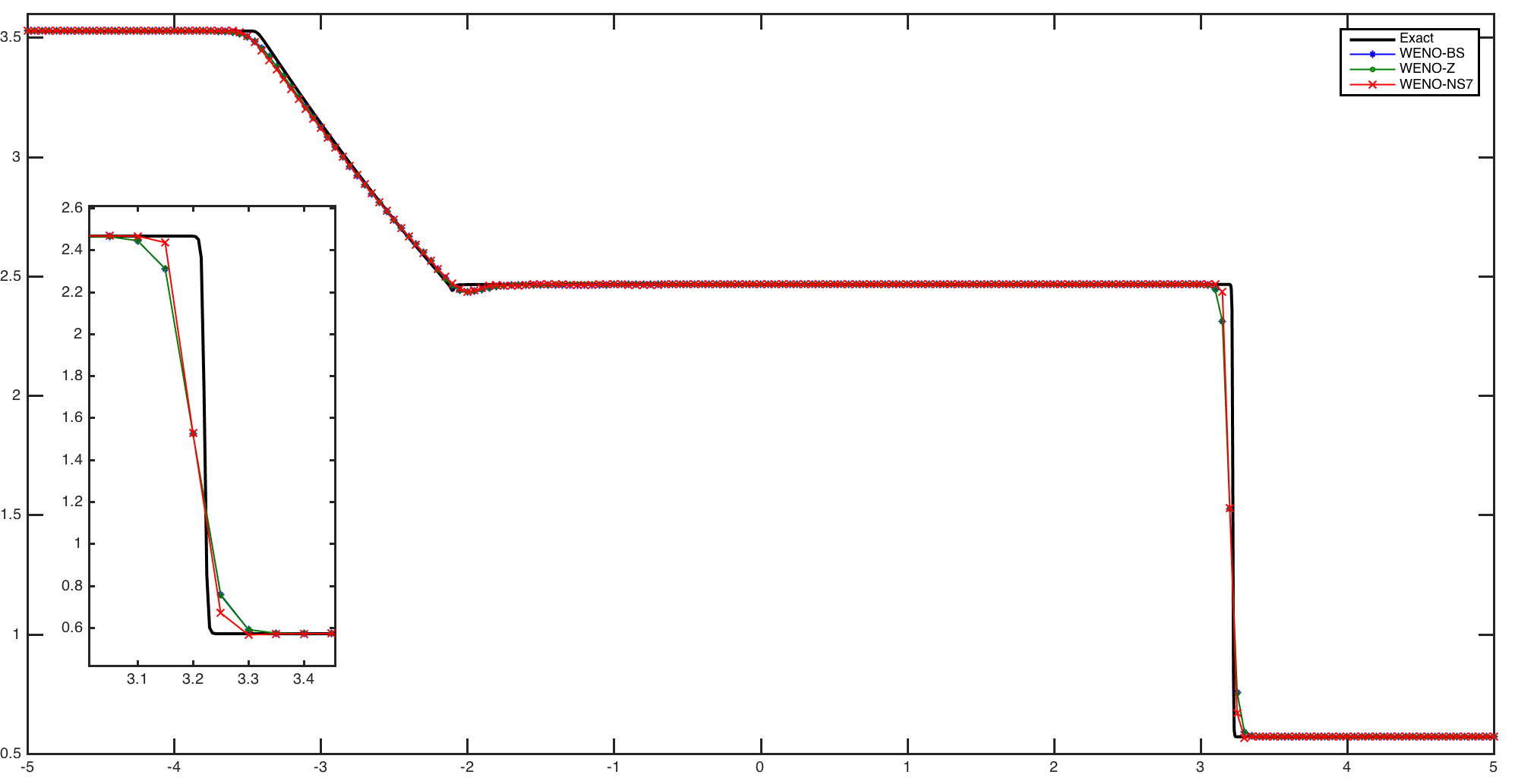}\caption{\label{fig:Pres-Lax}Pressure profile for Lax initial condition.}
\end{figure}

\begin{figure}[H]
\centering{}\includegraphics[width=14cm,height=6cm]{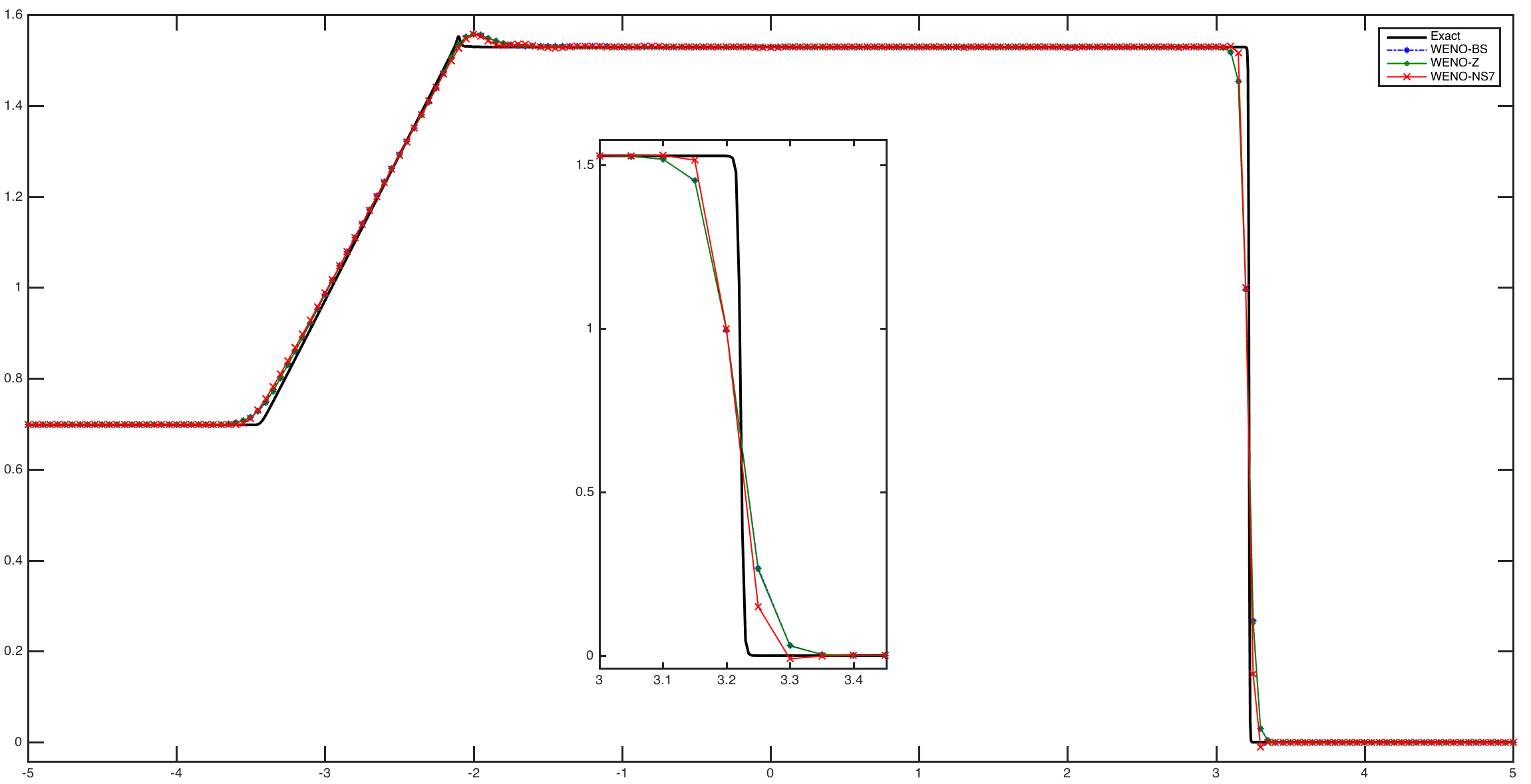}\caption{\label{fig:vel-Lax}Velocity profile for Lax initial condition.}
\end{figure}

The numerical solutions are computed on the computational domain $-5\leq x\leq5$
with $200$ grid cells up to time $t=1.3$ . Figures \ref{fig:Den-Lax},
\ref{fig:Pres-Lax} and \ref{fig:vel-Lax} displays that the density,
pressure and velocity profiles against the reference solution. WENO-NS7
performs better in this case too.

\subsubsection{1D shock entropy wave interaction problem}
\textbf{Example 7:}
Consider the system (\ref{eq:60}) on the spatial domain $x\in[-5,5]$
, with the Riemann data
\begin{equation}
(\rho,u,p)=\begin{cases}
(3.857143,2.629369,10.33333), & \text{ if }-5\leq x<-4,\\
(1+\epsilon sin(kx),0.000,1.000), & \text{ if }-4\leq x\leq5.
\end{cases}\label{60}
\end{equation}
where $\epsilon$ and $k$ are the amplitude and wave number of the
entropy wave respectively, chosen as $\epsilon=0.2$ and $k=5$ .
This problem is known as shock entropy wave interaction problem \cite{woodward},
the solution has a right-moving supersonic (Mach 3) shock wave which
interacts with sine waves in a density disturbance that generates
a flow field with both smooth structures and discontinuities. This
flow induces wave trails behind a right-going shock with wave numbers
higher than the initial density-variation wave number $k$. The initial
condition contains a jump discontinuity at $x=-4,$ especially the
initial density profile has oscillations on $[-4,5]$.

\begin{figure}[H]
\centering{}\includegraphics[width=14cm,height=8cm]{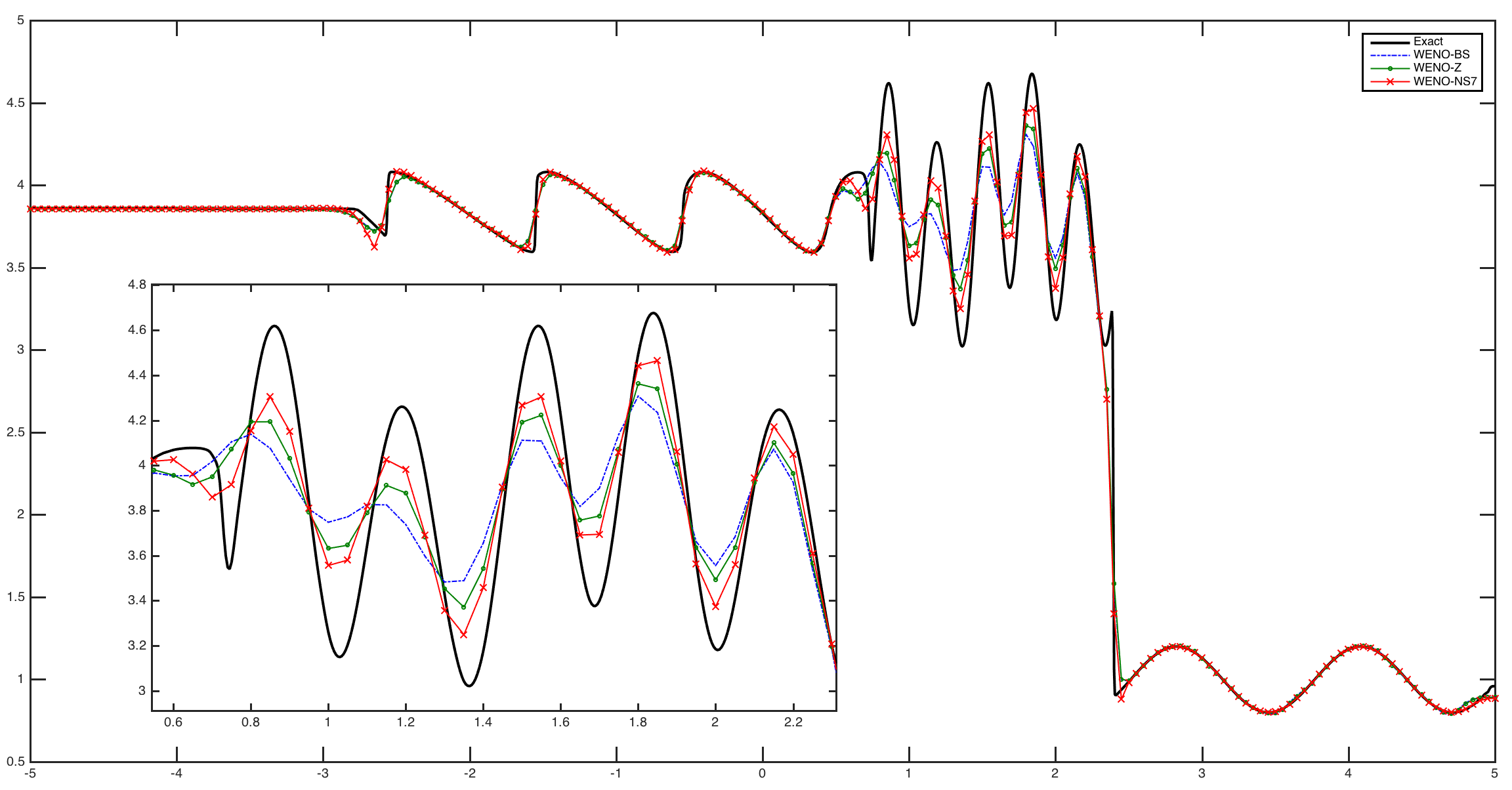}\caption{\label{fig:Shock-entropy-200}Shock entropy wave interaction test
with 200 grid points}
\end{figure}

\begin{figure}[H]
\centering{}\includegraphics[width=14cm,height=7cm]{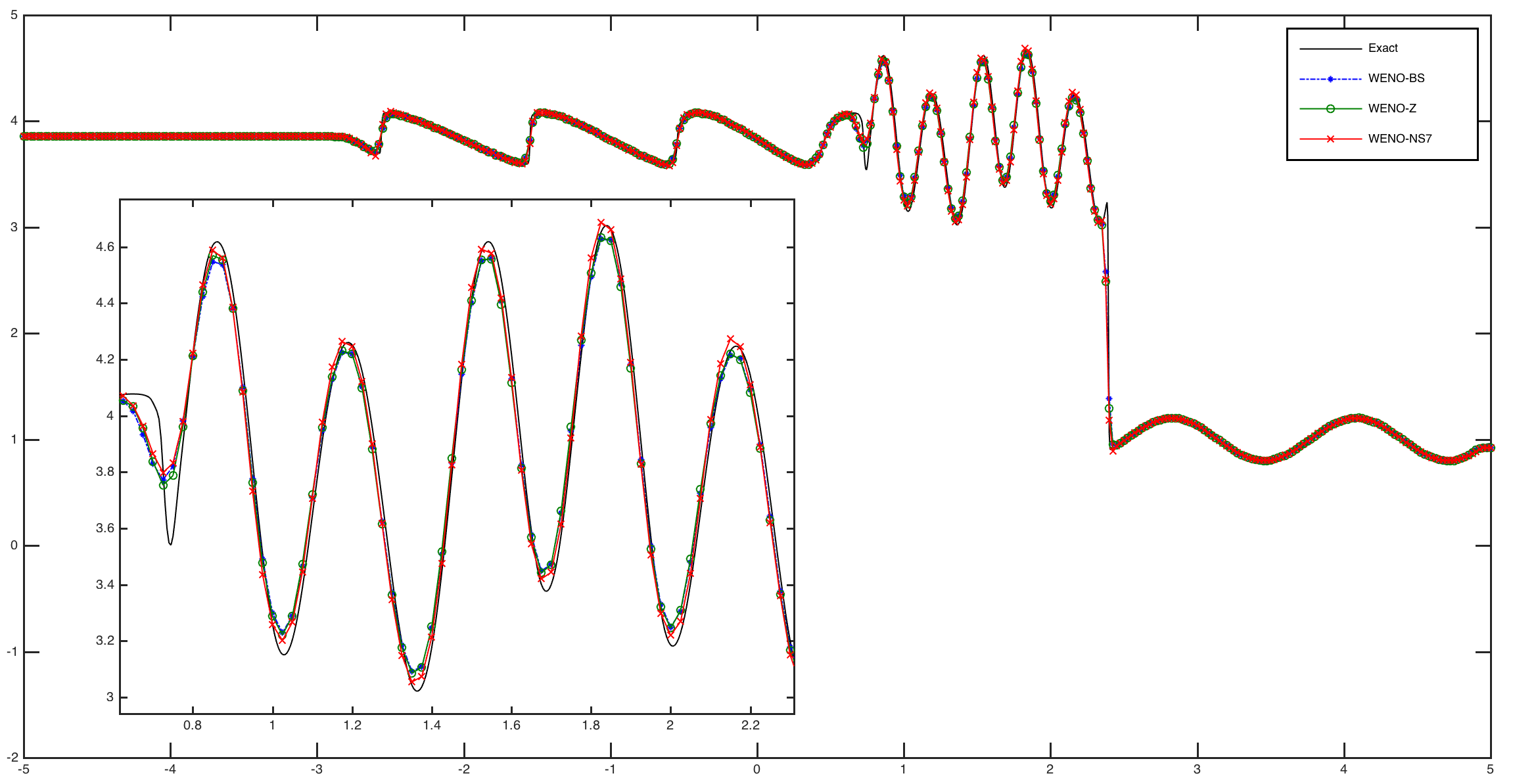}\caption{\label{fig:Shock-entropy-400}Shock entropy wave interaction test
with 400 points}
\end{figure}

The numerical results are computed with $200$ and $400$ spatial
grid points up to time $t=1.8,$ the results are plotted against the
reference solution in figures \ref{fig:Shock-entropy-200} and \ref{fig:Shock-entropy-400}
respectively. The proposed scheme WENO-NS7 performs better than other
seventh order schemes, which can be observed from these figures.

\subsection{Two-Dimensional Euler equations}

\hspace{0.6 cm}In this section, we apply the proposed scheme to two-dimensional
problem in cartesian coordinates. The governing two-dimensional compressible
Euler equations are given by
\begin{equation}
U_{t}+F(U)_{x}+G(U)_{y}=0,\label{eq:2D}
\end{equation}
where
\[
U=(\rho,\rho u,\rho v,E)^{T},\,F(U)=(\rho u,P+\rho u^{2},\rho uv,u(E+P))^{T},G(U)=(\rho v,\rho uv,P+\rho v^{2},v(E+P))^{T}.
\]
Here $\rho,u,v$ are density, $x$-velocity component and $y$-velocity
component respectively.

\begin{figure}[H]
\centering{}\includegraphics[width=16cm,height=6cm]{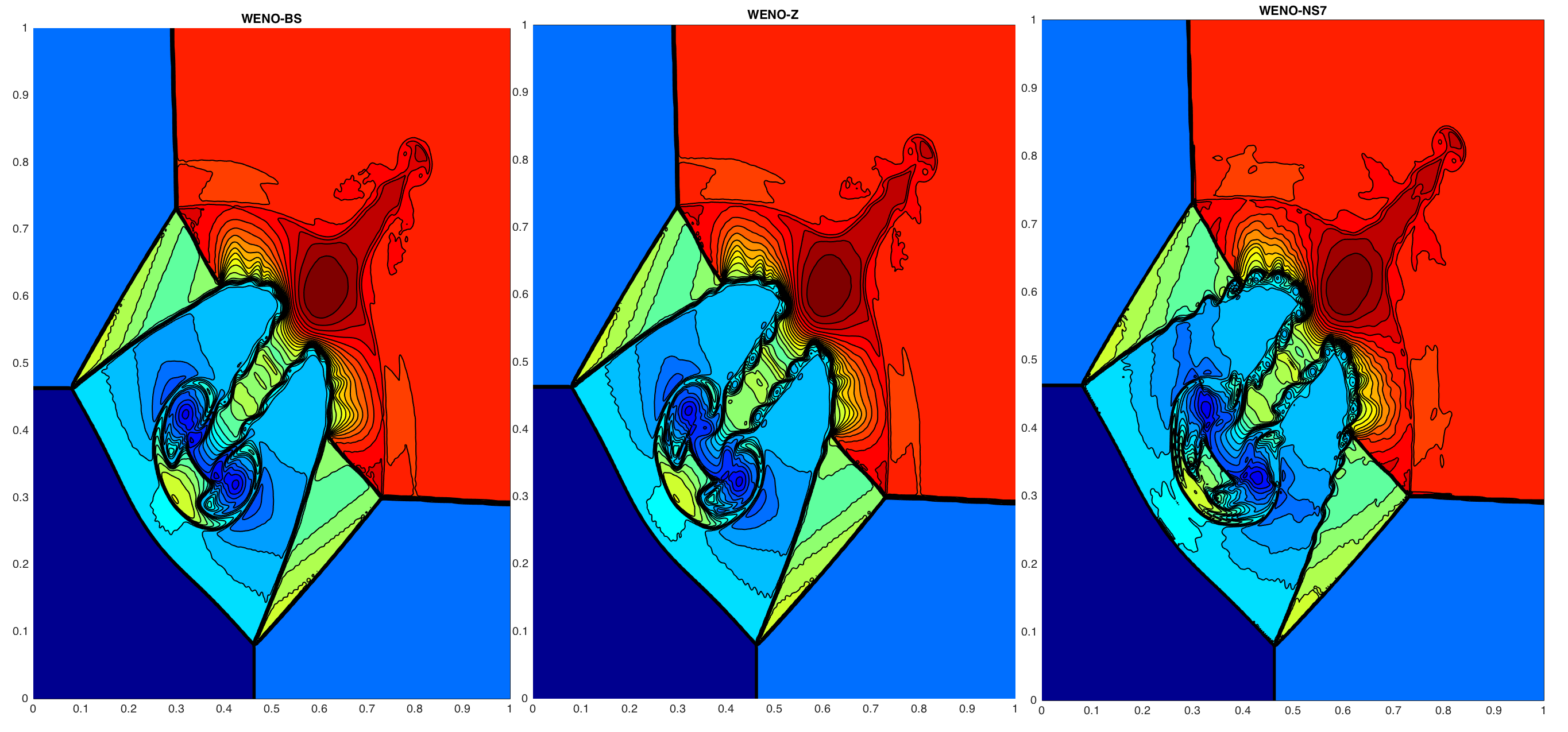}\caption{\label{fig:2d density}Density profiles of 2D problem with the mesh
$\Delta x=\Delta y=1/400$}
\end{figure}

The total energy $E$ can be obtained from the equation

\[
p=(\gamma-1)(E-\frac{1}{2}\rho(u^{2}+v^{2})),
\]
where $p$ is the pressure and $\gamma$ is the ratio of specific
heats.\\
\textbf{Example 8:}
For the PDE (\ref{eq:2D}) on the domain $[0,1]\times[0,1]$, consider
the initial condition
\[
(\rho,u,v,p)=\begin{cases}
(1.5,0,0,1.5) & \text{ if }0.8\leq x\leq1,0.8\leq y\leq1,\\
(0.5323,1.206,0,0.3) & \text{ if }0\leq x<0.8,0.8\leq y\leq1,\\
(0.138,1.206,1.206,0.029) & \text{ if }0\leq x<0.8,0\leq y<0.8,\\
(0.5323,0,1.206,0.3) & \text{ if }0.8<x\leq1,0\leq y\leq0.8,
\end{cases}
\]
given by Rinne et. al. in \cite{Rinne}, here each quadrant is divided
by lines $x=0.8$ and $y=0.8$. This initial condition produce a narrow
jet by the merge of four shocks.
The numerical solution is computed on the mesh of $400\times400$
grid points up to time $t=0.8$ with Dirichlet boundary conditions
and is displayed in figure \ref{fig:2d density}. An examination of
this result reveal that WENO-NS7 yields a better solution of the complex
structure when compares to WENO-BS and WENO-Z schemes.

\section{Conclusions}

\hspace{0.6cm} In this article, a new seventh-order weighted essentially
non-oscillatory scheme is presented to approximate the solution of
nonlinear hyperbolic conservation laws, named as WENO-NS7 scheme.
Based on undivided differences of derivatives of an interpolation
polynomials over a stencil, the local as well as global smoothness
indicators are constructed, measured in $L_{1}-$sense, to approximate
the nonlinear weights. We analyzed the scheme in the presence of critical
points and found that the proposed scheme possesses the desirable
order of accuracy when the first and second-order derivative vanishes.
The approximate solutions of the one and two-dimensional hyperbolic
conservation laws are simulated with the proposed WENO-NS7 scheme
and compared it with WENO-BS and WENO-Z schemes. Numerical experiments
show that the proposed WENO-NS7 scheme yields better approximation
in comparison to the WENO-BS and WENO-Z schemes. As a future work,
working on the construction of higher order schemes (order greater
than seven).

\section{Appendix}

For completeness, here we are describing briefly the WENO-BS and WENO-Z
schemes, these schemes are used in numerical comparison.

\subsection{WENO-BS}

\hspace{0.6 cm}For the reconstruction of flux (\ref{eq:recon}),
Balsara and Shu \cite{balsara and shu8} has taken the non-linear
weights $\omega_{k}'s$ as

\begin{equation}
\omega_{k}=\frac{\alpha_{k}}{\sum_{q=0}^{3}\alpha_{q}},\alpha_{k}=\frac{d_{k}}{(\epsilon+\beta_{k})^{p}},k=0,1,2,3,
\end{equation}
where $0<\epsilon\ll1$ is introduced to prevent $\alpha_{k}\rightarrow\infty,$
the smoothness indicator $\beta_{k},$ measures the smoothness of
a solution over a particular stencil and $p\geq1$ is a parameter.
The value of $p$ is chosen such that the nonlinear weights, in non-smooth
regions approaches to zero at an accelerated rate as $\Delta x\rightarrow0$.
The suggested smoothness indicators $\beta_{k}$ are the same as that
of WENO-JS \cite{Jiang and shu7}, which are given by
\begin{equation}
\beta_{k}=\sum_{q=1}^{3}\triangle x^{2q-1}\intop_{x-\frac{\triangle x}{2}}^{x+\frac{\triangle x}{2}}\left(\frac{d^{q}\hat{f}^{k}}{dx^{q}}\right)^{2}dx.\label{eq:bks}
\end{equation}
The explicit form of these smoothness indicators are as follows:

\begin{eqnarray*}
\begin{split}\beta_{0}= & f_{j-3}(547f_{j-2}-3882f_{j-2}+4642f_{j-1}-1854f_{j})+f_{j-2}(7043f_{j-2}-17246f_{j-1}+7042f_{j})+\\
 & f_{j-1}(11003f_{j-1}-9402f_{j})+2107f_{j}^{2},\quad
\end{split}
\end{eqnarray*}
\begin{eqnarray*}
\begin{split}\beta_{1}= & f_{j-2}(267f_{j-2}-1642f_{j-1}+1602f_{j}-494f_{j+1})+f_{j-1}(2843f_{j-1}-5966f_{j}+1922f_{j+1})+\\
 & f_{j}(3443f_{j}-2522f_{j+1})+547f_{j+1}^{2},\quad
\end{split}
\end{eqnarray*}
\begin{eqnarray*}
\begin{split}\beta_{2}= & f_{j-1}(547f_{j-1}-2522f_{j}+1922f_{j+1}-494f_{j+2})+f_{j}(3443f_{j}-5966f_{j+1}+1602f_{j+2})+\quad\\
 & f_{j+1}(2843f_{j+1}-1642f_{j+2})+267f_{j+2}^{2},\qquad\quad
\end{split}
\end{eqnarray*}
\begin{eqnarray*}
\begin{split}\beta_{3}= & f_{j}(2107f_{j}-9402f_{j+1}+7042f_{j+2}-1854f_{j+3})+f_{j+1}(11003f_{j+1}-17246f_{j+2}+4642f_{j+3})+\\
 & f_{j+2}(7043f_{j+2}-3882f_{j+3})+547f_{j+3}^{2}.
\end{split}
\end{eqnarray*}

\subsection{WENO-Z }

\hspace{0.6 cm}Borges et al. \cite{Borges caramona10} redefined
the non-linear weights of the WENO-JS scheme \cite{Jiang and shu7}
by introducing a global smoothness indicator $\tau,$ as
\begin{equation}
\omega_{k}^{z}=\frac{\alpha_{k}^{z}}{\sum_{q=0}^{3}\alpha_{q}^{z}},\,\alpha_{k}^{z}=d_{k}\bigg(1+\bigg[\frac{\tau}{\beta_{k}+\epsilon}\bigg]^{p}\bigg).\label{23}
\end{equation}
The idea here is to get the nonlinear weights $\omega_{k}'s$ close
to the ideal weights $d_{k}'s,$ the scheme is known as WENO-Z scheme.

For the seventh order WENO-Z scheme, Castro et. al. \cite{Castro et al}
defined the global smoothness indicator as
\begin{equation}
\tau_{7}=|\beta_{0}-\beta_{1}-\beta_{2}+\beta_{3}|.\label{eq:25}
\end{equation}
whose truncation error is $O\left(\Delta x^{7}\right)$ and in \cite{W S DON }
the author's defined another global smoothness indicator
\begin{equation}
\tau_{7}=|\beta_{0}+3\beta_{1}-3\beta_{2}-\beta_{3}|\label{eq:STD}
\end{equation}
whose truncation error is $O\left(\Delta x^{8}\right).$ The numerical
convergence of the WENO-Z scheme, when used with the high-order global
smoothness indicator (\ref{eq:STD}) achieves better order of accuracy
in comparison to the use of the global smoothness indicator (\ref{eq:25}).
For the numerical comparison, we've used the global smoothness indicator
(\ref{eq:STD}) with the nonlinear weights (\ref{23}) to compute
the WENO-Z scheme. The smoothness indicators $\beta_{k}'s$ in (\ref{eq:25})
and (\ref{eq:STD}) are as in (\ref{eq:bks}).


\begin{thebibliography}{}
\bibitem[1]{Arandiga} Arandiga, F., Baeza, A.,  Belda, A.M., Mulet, P.:
Analysis of WENO schemes for full and global accuracy. SIAM J. Numer.
Anal. \textbf{49}, 893-915(2011).

\bibitem[2]{balsara and shu8} Balsara, D.S.,  Shu, C.W.:  Monotonicity
preserving weighted essentially non-oscillatory schemes with increasingly
high order of accuracy. J. Comput. Phys. \textbf{160}, 405-452(2000).

\bibitem[3]{Borges caramona10} Borges, R., Carmona, M., Costa, B., Don,  W.S.: An improved weighted essentially non-oscillatory scheme for hyperbolic
conservation laws. J. Comput. Phys. \textbf{227}, 3191-3211(2008).

\bibitem[4]{Castro et al}Castro, M., Costa,   B., Don, W.S.: High order
weighted essentially nonoscillatory WENO-Z schemes for hyperbolic
conservation laws. J. Comput. Phys. \textbf{230}, 766-792(2011).

\bibitem[5]{W S DON }Don, W.S., Borges,  R.:  Accuracy of the weighted
essentially non-oscillatory conservative finite difference schemes.
J. Comput. Phys. \textbf{250}, 347-372(2013).

\bibitem[6]{P fan et.al.14}Fan, P., Shen,  Y., Tian,  B., Yang, C.:  A new
smoothness indicator for improving the weighted essentially non-oscillatory
scheme. J. Comput. Phys. \textbf{269}, 329-354(2014).

\bibitem[7]{P fan15} Fan, P.:  High order weighted essentially non
oscillatory WENO-schemes for hyperbolic conservation laws. J. Comput.
Phys. \textbf{269}, 355-285(2014).

\bibitem[8]{Gerolymus}Gerolymos, G.A., Senechal, D., Vallet, I.:  Very
high order WENO schemes. J. Comput. Phys. \textbf{228}, 8481-8524(2009).

\bibitem[9]{Godunov}Godunov, S.K.:  A finite-difference scheme for
the numerical computation of discontinuous solutions of the equations
of fluid dynamics. Matthematicheskii sbornik, \textbf{47}, 271-290(1959).

\bibitem[10]{Gottilieb}Gottlieb, S.:  On high order strong stability
preserving Runge\textendash Kutta and multi- step time discretizations.
J. Sci. Comput. \textbf{25}, 105-28(2005) .

\bibitem[11]{Ha et.al13}Ha, Y., Kim, C.H., Lee,  Y.J., Yoon,  J.:  An improved
weighted essentially non-oscillatory scheme with a new smoothness
indicator. J. Comput. Phys. \textbf{232}, 68-86(2013).

\bibitem[12]{Ha83}Harten, A.:  High resolution schemes for hyperbolic
conservation laws. J. Comput. Phys. \textbf{49}, 357-393(1983).

\bibitem[13]{Ha84}Harten, A.:  On a class of high resolution total-variation-stable
finite-difference schemes. SIAM J. Numer. Anal. \textbf{21},1(1984).

\bibitem[14]{HOEC86}Harten, A., Osher, S.,  Engquist, B., Chakravarthy, S.R.:
Some results on uniformly high-order accurate essentially nonoscillatory
schemes. App. Numer. Math.\textbf{2}, 347-377(1986).

\bibitem[15]{HO87}Harten, A., Osher, S.:  Uniformly high-order accurate
non oscillatory schemes, I. SIAM J. Numer. Anal. \textbf{24}, 279-309(1987).

\bibitem[16]{HO97}Harten, A., Engquist, B., Osher,  S., Chakravarthy,  S.R.:
Uniformly high order accurate non-oscillatory schemes, III. J. Comput.
Phys. \textbf{131}, 3-47(1997).

\bibitem[17]{henrick aslam powers5}Henrick, A.K., Aslam,  T.D.,
Powers, J.M.: Mapped weighted essentially non-oscillatory schemes: Achieving
optimal order near critical points. J. Comput. Phys. \textbf{207}, 542-567(2005).

\bibitem[18]{MWENO-Z}Hu, F., Wang,  R., Chen,  X.:  A modified fifth-order
WENOZ method for hyperbolic conservation laws. J. Comput. and App.
math. \textbf{303}, 56-68(2016).

\bibitem[19]{Jiang and shu7}Jiang,  G.S., Shu,  C.W.:  Efficient implementation
of Weighted ENO schemes. J. Comput. Phys. \textbf{126}, 202-228(1996).

\bibitem[20]{KIM et al}Kim, C.H.,   Ha, Y., Yoon, J.:  Modified Non-linear
Weights for Fifth-Order Weighted Essentially Non-oscillatory Schemes.
J. Sci. Comput. \textbf{67}, 299-323(2016).

\bibitem[21]{PDLax}Lax, P.D.: Weak solutions of nonlinear hyperbolic equations and their numerical computation. Commun. Pure Appl. Math. \textbf{7}, 159-193(1954).

\bibitem[22]{XD Liu}Liu, X.D., Osher,  S., Chan,  T.: Weighted Essentially
non-oscillatory schemes. J. Comput. Phys. \textbf{115}, 200-212(1994).

\bibitem[23]{OC84}Osher, S., Chakravarthy,  S.R.:  High resolution schemes
and the entropy condition. SIAM J. Numer. Anal. \textbf{21},5(1984).

\bibitem[24]{RathanRaju}Rathan, S., Naga Raju,  G.:  A modified fifth-order
WENO scheme for  hyperbolic conservation laws. arxiv preprint.

\bibitem[25]{Rinne} Schulz-Rinne, C.W., Collins, J.P., Glaz,  H.M.: Numerical
solution of the riemann problem for two-dimensional gas dynamics.
SIAM J. Sci. Compu. \textbf{14}, 1394-1414(1993).

\bibitem[26]{serna}Serna, S., Marquina,  A.:  Power-ENO methods: a fifth-order
accurate weighted power ENO method. J. Comput. Phys. \textbf{194}, 632-658(2004).

\bibitem[27]{ShenZha}Shen, Y.Q., Zha,  G.C.:  A Robust Seventh-order
WENO Scheme and Its Applications. AIAA-2008-0757(2008).

\bibitem[28]{ShenZha2}Shen, Y.Q.,  Zha,  G.C.: Improved Seventh-Order
WENO Scheme. AIAA-2010-1451(2010).

\bibitem[29]{shu notes}Shu, C.W.:  Essentially non-oscillatory and
weighted essentially non-oscillatory schemes for hyperbolic conservation
laws. In Advanced numerical approximation of nonlinear hyperbolic
equations, Lecture Notes in Mathematics,Berlin, Springer-Verlag. \textbf{1697}, 325-432 (1998).


\bibitem[30]{c w shu16}Shu, C.W.:  High order weighted essentially
non oscillatory schemes for convection dominated problems. SIAM Review.
\textbf{51},82-126(2009).

\bibitem[31]{osherSHU}Shu, C.W., Osher, S.:  Efficient implementation
of essentially non-oscillatory shock-capturing schemes. J. Comput.
Phys. \textbf{77}, 439-471(1988).

\bibitem[32]{Shu-osher1}Shu, C.W., Osher,  S.:  Efficient implementation
of essentially non-oscillatory shock-capturing schemes II. J. Comput.
Phys. \textbf{83}, 32-78(1989).

\bibitem[33]{Sod}Sod, G.A.:  A survey of several finite difference
methods for systems of nonlinear hyperbolic conservation laws. J.
Comput. Phys.\textbf{107}, 1-31(1978).

\bibitem[34]{woodward}Woodward, P.,  Colella,  P.: The numerical simulation
of two-dimensional fluid flow with strong shocks. J. Comput. Phys.
\textbf{54}, 115-173(1984).


\end{thebibliography}
\end{document}